%% file: article.tex
\documentclass[final]{siamltex}

\usepackage{xspace}
\usepackage{multirow}
\usepackage{array}
\usepackage{ifthen}
\usepackage{amsmath}
\usepackage{amsfonts}
\usepackage{color}
\usepackage{verbatim}
\usepackage{relsize}
\usepackage{mathdots}
\usepackage{wrapfig}
\usepackage{array}
\usepackage{rotating}
\usepackage{url}
\usepackage{tikz}
\usepackage{pgfplots}
\pgfplotsset{compat=newest}

\newcommand{\href}[1]{\url{#1}}

\input{MOOC_macros}
\newtheorem{remark}{Remark}

\newcommand{\pgm}[1]{{\color{black} #1}}

\newcommand{\RvdG}[1]{{\color{black} #1}}
\title{
Householder QR Factorization With
Randomization for Column Pivoting (HQRRP)
}

\author{
Per-Gunnar Martinsson%
\footnotemark[1]
\and
Gregorio Quintana Ort\'{\i}%
\footnotemark[2]
\and
Nathan Heavner\footnotemark[1]
\and
Robert van de Geijn%
\footnotemark[3]
}

\date{Draft \\ \today}
\begin{document}

\maketitle

\renewcommand{\thefootnote}{\fnsymbol{footnote}}
\footnotetext[1]{Department of Applied Mathematics, University of Colorado at Boulder, 526 UCB, Boulder, CO 80309-0526, USA
}
\footnotetext[2]{Depto.~de Ingenier\'{\i}a y Ciencia de Computadores, Universitat Jaume I, 12.071 Castell\'on, Spain
}
\footnotetext[3]{
Department of Computer Science and  Institute for Computational Engineering and Sciences, The University of Texas at Austin,
Austin, TX.
%\newline Email: rvdg@cs.utexas.edu.
}

\begin{abstract}
\input{00abstract}
\end{abstract}

\input{body}

\subsection*{Acknowledgments}

The research reported was supported by DARPA,
under the contract N66001-13-1-4050, and by the NSF, under the contracts DMS-1407340, DMS-1620472,
and ACI-1148125/1340293.
\textit{Any opinions, findings and conclusions or recommendations
 expressed in this material are those of the author(s) and do not
 necessarily reflect the views of the National Science Foundation
 (NSF).}

\bibliographystyle{plain}
\bibliography{biblio,main_bib}

\end{document}

%% file: MOOC_macros.tex
\input flatex_color

\newcommand{\Answer}[1]{%
\ifthenelse{\boolean{ShowAnswers}}{%
{\color{blue}%
\vspace{0.05in}%
\noindent%
{\bf Answer:} #1}}
{}%
}

\newcommand{\PartAnswer}[1]{%
\ifthenelse{\boolean{ShowAnswers}}{%
{\color{blue} #1}}%
{}
}

\newcommand{\QuestionAnswer}[2]{%
\ifthenelse{\boolean{ShowAnswers}}{%
{\color{blue} #2}}%
{\color{black} #1}
}

\newcommand{\ShowAnswer}[1]{
\ifthenelse{\boolean{ShowAnswers}}{
{\color{blue} #1}}
{\phantom{#1}}
}

{
\vspace{0.1in}
\end{minipage}
\\ \whline
\end{tabular}
}

\newenvironment{unboxedexercise}{
\addtocounter{homeworkcounter}{1}%
\noindent%
\begin{exercise}
\rm 
}
{
\end{exercise}%
}

\newcolumntype{I}{!{\vrule width 1.5pt}}
\newlength\savedwidth
\newcommand\whline{\noalign{\global\savedwidth\arrayrulewidth
                            \global\arrayrulewidth 1.5pt}%
           \hline
           \noalign{\global\arrayrulewidth\savedwidth}}

\newcommand{\NoShow}[1]{}

\newcommand{\becomes}{:=}

\newcommand{\Cn}{\mathbb{C}^n}

\newcommand{\Cmxn}{\mathbb{C}^{m \times n}}

\newcommand{\housev}{{\rm Housev}}

\newcommand{\gemm}{\mbox{\sc gemm}}
\newcommand{\HQR}{HQR}
\newcommand{\HQRP}{HQRP}
\newcommand{\HQRRP}{HQRRP}
\newcommand{\QRP}{QRP}

%% file: flatex_color.tex
\usepackage{ifthen}
\usepackage{colortbl}
\usepackage{array}
\usepackage{color}
\usepackage{soul}
\usepackage{amssymb}

\newboolean{IsWide}
\setboolean{IsWide}{true}

\newcommand{\FlaTwoByTwo}[4]{
\left(
\begin{array}{c I c}
#1 & #2 \\ \whline
#3 & #4
\end{array}
\right)
}

\newcommand{\FlaTwoByTwoSingleLine}[4]{
\left(
\begin{array}{c | c}
#1 & #2 \\ \hline
#3 & #4
\end{array}
\right)
}

\newcommand{\FlaTwoByOne}[2]{
\left(
\begin{array}{c}
#1 \\ \whline
#2
\end{array}
\right)
}

\newcommand{\FlaTwoByOneSingleLine}[2]{
\left(
\begin{array}{c}
#1 \\ \hline
#2
\end{array}
\right)
}

\newcommand{\FlaOneByTwo}[2]{
\left(
\begin{array}{c I c}
#1 & #2
\end{array}
\right)
}

\newcommand{\FlaOneByTwoSingleLine}[2]{
\left(
\begin{array}{c | c}
#1 & #2
\end{array}
\right)
}

\newcommand{\FlaThreeByThreeTL}[9]{
\left(
\begin{array}{c | c I c}
#1 & #2 & #3 \\ \hline
#4 & #5 & #6 \\ \whline
#7 & #8 & #9
\end{array}
\right)
}

\newcommand{\FlaThreeByThreeBR}[9]{
\left(
\begin{array}{c I c | c}
#1 & #2 & #3 \\ \whline
#4 & #5 & #6 \\ \hline
#7 & #8 & #9
\end{array}
\right)
}

\newcommand{\FlaOneByThreeR}[3]{
\left(
\begin{array}{c I c | c}
#1 & #2 & #3
\end{array}
\right)
}

\newcommand{\FlaOneByThreeL}[3]{
\left(
\begin{array}{c | c I c}
#1 & #2 & #3
\end{array}
\right)
}

\newcommand{\FlaThreeByOneT}[3]{
\left(
\begin{array}{c}
#1 \\ \hline
#2 \\ \whline
#3
\end{array}
\right)
}

\newcommand{\FlaThreeByOneB}[3]{
\left(
\begin{array}{c}
#1 \\ \whline
#2 \\ \hline
#3
\end{array}
\right)
}

\newcommand{\FlaPartition}[2]{
\ifthenelse{\boolean{IsWide}}{{\bf partition } \hspace{-1em} #1 \hspace{-1em} #2}
{{\bf partition } \+ \\ #1 \+ \\ #2 \- \-}
}

\newcommand{\FlaRepartition}[2]{
\ifthenelse{\boolean{IsWide}}{{\bf repartition } \hspace{-1em} #1 \hspace{-1em} #2}
{{\bf repartition } \+ \\ #1 \+ \\ #2 \- \-}
}

\newcommand{\FlaStartCompute}{%
\setlength{\unitlength}{1.6in}%
\begin{picture}(3,0.01)
\put(0,0){\line(1,0){3}}
\put(0,0.01){\line(1,0){3}}
\end{picture}%
}

\newcommand{\FlaEndCompute}{%
\noindent%
\setlength{\unitlength}{1.6in}%
\begin{picture}(3,0.01)
\put(0,0){\line(1,0){3}}
\put(0,0.01){\line(1,0){3}}
\end{picture}%
}

\newcommand{\FlaStartComputeShorter}{
{\color{red}
\setlength{\unitlength}{1.6in}
\begin{picture}(1.2,0.01)
\put(0,0){\line(1,0){1.2}}
\put(0,0.01){\line(1,0){1.2}}
\end{picture}
}
}

\newcommand{\FlaEndComputeShorter}{
{\color{red}
\setlength{\unitlength}{1.6in}
\begin{picture}(1.2,0.01)
\put(0,0){\line(1,0){1.2}}
\put(0,0.01){\line(1,0){1.2}}
\end{picture}
}
}

\newcommand{\operation}{ [ D, E, F, \ldots ] \becomes {\rm op}( A, B, C, D, \ldots ) }
\newcommand{\routinename}{ [ D, E, F, \ldots ] \becomes {\rm op}( A, B, C, D, \ldots ) }
\newcommand{\routinecost}{ X }
\newcommand{\precondition}{ Q }
\newcommand{\postcondition}{ R }
\newcommand{\invariant}{ P }
\newcommand{\costinv}{ \  }

\newcommand{\guard}{ R }
\newcommand{\partitionings}{
\begin{minipage}{3in}
$ S_I $
\end{minipage}
}
\newcommand{\initialize}{}

\newcommand{\partitionsizes}{ \hspace{ 3.25in} }

\newcommand{\blocksize}{1}

\newcommand{\repartitionings}{
\begin{minipage}[t]{3in}
\ \\
\ \\
\ \\
\end{minipage}
}

\newcommand{\repartitionsizes}{ \hspace{ 3.25in} }
\newcommand{\moveboundaries}{
\begin{minipage}[t]{3in}
\ \\
\ \\
\ \\
\end{minipage}
}
\newcommand{\beforeupdate}{
$ \QBefore $
}
\newcommand{\afterupdate}{
$ \QAfter $
}
\newcommand{\update}{%
\begin{minipage}[t]{4in}
$ S_U $
\end{minipage}
}

\newcommand{\resetsteps}{
\renewcommand{\blocksize}{1}
\renewcommand{\operation}{ [ D, E, F, \ldots ] \becomes {\rm op}( A, B, C, D, \ldots ) }
\renewcommand{\routinename}{ [ D, E, F, \ldots ] \becomes {\rm op}( A, B, C, D, \ldots ) }
\renewcommand{\routinecost}{ 0 }
\renewcommand{\precondition}{ \PPre }
\renewcommand{\postcondition}{ \PPost }
\renewcommand{\invariant}{ \PInv }
\renewcommand{\costinv}{ \  }
\renewcommand{\guard}{ G }
\renewcommand{\partitionings}{ %
\begin{minipage}[t]{3in}
\ \\
\end{minipage}
}
\renewcommand{\partitionsizes}{ \hspace{ 3.25in} }
\renewcommand{\repartitionings}{%
\begin{minipage}[t]{3in}
\ \\
\end{minipage}
}
\renewcommand{\repartitionsizes}{ \hspace{ 3.25in} }
\renewcommand{\moveboundaries}{%
\begin{minipage}[t]{3in}
\ \\
\end{minipage}
}
\renewcommand{\beforeupdate}{
\QBefore
}
\renewcommand{\afterupdate}{
\QAfter
}
\renewcommand{\update}{
$ S_U $
}
}

\newcommand{\WSguard}{
$ \guard $
}

\newcommand{\WSpartition}{%
\begin{minipage}[t]{3.0in}%
\begin{tabbing}
ind \= ind \= \kill
{\bf \color{blue} Partition}
\partitionings \+ \\
{\bf \color{blue} where } \hspace*{-2ex} \partitionsizes
\end{tabbing}
\end{minipage}
}

\newcommand{\WSpartitionNarrow}{%
\begin{minipage}[t]{2.25in}%
\begin{tabbing}
id \= id \= \kill
{\bf \color{blue} Partition}
\partitionings \+ \\
{\bf \color{blue} where }
\begin{minipage}[t]{2.5in}
\partitionsizes
\end{minipage}
\end{tabbing}
\end{minipage}
}

\newcommand{\WSinitialize}{
\initialize
}

\newcommand{\WSrepartition}{
\begin{minipage}[t]{3in}
\ifthenelse{ \equal{\blocksize}{1} }{}
{%
\ifthenelse{ \equal{\blocksize}{blank} }{~}
{{\bf \color{blue} Determine block size} $ \blocksize $} \\
}
{\bf \color{blue} Repartition}
\begin{tabbing}
in \= in \= \+ \kill
\repartitionings \+ \\
{\bf \color{blue} where } \hspace*{-2ex} \repartitionsizes
\end{tabbing}
\end{minipage}
}

\newcommand{\WSrepartitionNarrow}{
\begin{minipage}[t]{2.05in}
\ifthenelse{ \equal{\blocksize}{1} }{\phantom{Determine}~}
{%
\ifthenelse{ \equal{\blocksize}{blank} }{}
{{\bf \color{blue} Determine block size} $ \blocksize $} \\
}
{\bf \color{blue} Repartition}
\begin{tabbing}
i \= i \= \+ \kill
\repartitionings 
\+ \\
{\bf \color{blue} where }
\begin{minipage}[t]{1.5in}
\repartitionsizes
\end{minipage}
\end{tabbing}
\end{minipage}
}

\newcommand{\WSmoveboundary}{%
\begin{minipage}[t]{4in}%
{\bf \color{blue} Continue with}
\begin{tabbing}
ind \= \+ \kill
\moveboundaries
\end{tabbing}
\end{minipage}
}

\newcommand{\WSmoveboundaryNarrow}{%
\begin{minipage}[t]{2.05in}%
{\bf \color{blue} Continue with}
\begin{tabbing}
i \= \+ \kill
\moveboundaries
\end{tabbing}
\end{minipage}
}

\newcommand{\WSupdate}{
\update
}

\newcommand{\FlaAlgorithmWithInit}{
\begin{center}
\begin{tabular}{| p{0.98\textwidth} |} \hline
{\bf \color{blue} Algorithm:} $\routinename$
\\ \whline
\WSinitialize \\
{\WSpartition} \\[0.3in]
{\bf \color{blue} while} \WSguard { \bf \color{blue} do} \\
\ \hspace{0.15in} \WSrepartition \\
\color{red} {\hspace{0.0in} \FlaStartCompute} \\
{\hspace{0.0in} \WSupdate} \\
\color{red} {\hspace{0.0in} \FlaEndCompute} \\
{\ \hspace{0.15in} \WSmoveboundary} \\
{{\bf \color{blue} endwhile}} \\ \hline
\end{tabular}
\end{center}
}

\newcommand{\FlaAlgorithmNarrow}{
\begin{center}
\begin{tabular}{| p{0.99\textwidth}|} \hline
{\bf \color{blue} Algorithm:} $\routinename$
\\ \whline
\color{black} \WSpartitionNarrow \\
{\bf \color{blue} while} \WSguard { \bf \color{blue} do} \\[-0.1in]
\ \hspace{0.0in} \color{black} \WSrepartitionNarrow \\[-0.1in]
{\hspace{0.0in} \FlaStartComputeShorter} \\
{\hspace{0.0in} \color{black} \WSupdate} \\[-0.1in]
{\hspace{0.0in} \FlaEndComputeShorter} \\[-0.1in]
{\ \hspace{0.0in} \color{black} \WSmoveboundaryNarrow} \\
{{\bf \color{blue} endwhile}} \\ \hline
\end{tabular}
\end{center}
}

\newcommand{\FlaAlgorithmWithInitNarrow}{
\begin{center}
\begin{tabular}{| p{0.99\textwidth}|} \hline
{\bf \color{blue} Algorithm:} $\routinename$
\\ \whline
\WSinitialize \\
\color{black} \WSpartitionNarrow \\
{\bf \color{blue} while} \WSguard { \bf \color{blue} do} \\[-0.1in]
\ \hspace{0.0in} \color{black} \WSrepartitionNarrow \\[-0.1in]
{\hspace{0.0in} \FlaStartComputeShorter} \\
{\hspace{0.0in} \color{black} \WSupdate} \\[-0.1in]
{\hspace{0.0in} \FlaEndComputeShorter} \\[-0.1in]
{\ \hspace{0.0in} \color{black} \WSmoveboundaryNarrow} \\
{{\bf \color{blue} endwhile}} \\ \hline
\end{tabular}
\end{center}
}

\newcommand{ \PPre }{ P_{\it pre} }
\newcommand{ \PPost }{ P_{\it post} }
\newcommand{ \PInv }{ P_{\it inv} }

\newcommand{ \QBefore }{ P_{\it before} }
\newcommand{ \QAfter }{ P_{\it after} }

%% file: 00abstract.tex
A fundamental problem when adding column pivoting to the Householder QR factorization is that only about half of the computation can be cast in terms of high performing matrix-matrix multiplications, which greatly limits the benefits that can be derived from so-called blocking of algorithms. This paper describes a technique for selecting groups of pivot vectors by means of randomized projections. It is demonstrated that the asymptotic flop count for the proposed method is $2mn^2 - (2/3)n^3$ for an $m\times n$ matrix, identical to that of the best classical unblocked Householder QR factorization algorithm (with or without pivoting). Experiments demonstrate acceleration in speed of close to an order of magnitude relative to the {\sc geqp3} function in LAPACK, when executed on a modern CPU with multiple cores. Further, experiments demonstrate that the quality of the randomized pivot selection strategy is roughly the same as that of classical column pivoting. The described algorithm is made available under Open Source license and can be used with LAPACK or libflame.
%The algorithm described
%has been implemented in the software package \texttt{HQRRP} (Householder QR with Randomization for
%Pivoting) that is publicly available through github, and includes a routine \texttt{dgeqp4}
%which is plug compatible with the LAPACK routine \texttt{dgeqp3}. 

%% file: body.tex
\newcommand{\mtx}[1]{#1}

\section{Introduction}

\input 01intro.tex

\section{Householder QR Factorization}

\input 02HQR

\section{Randomization to the Rescue}
\label{sec:HPQRRP}

\input 03HQRRP

\section{Experiments}
\label{sec:num}

\input 04experiments

\section{Conclusions and future work}

\input 05conclusion

\section{Software}

\input 06software

%% file: 01intro.tex
The QR factorization is a staple of linear algebra, with applications ranging from Linear Least-Squares solution of overdetermined systems to the identification of low rank approximation via the determination of an approximate orthonormal basis for the column space.  Standard algorithms for computing the QR factorization include Gram-Schmidt orthogonalization and those based on Householder transformations (reflectors).  When it is desirable for the QR factorization to also reveal the approximate rank of the original matrix, it is important that the elements of the diagonal of $ R $ be ordered with larger elements in magnitude appearing earlier.  In this case, column pivoting (swapping) is employed during the QR factorization, yielding QR factorization with column pivoting (\QRP).
It is well-known that the Householder QR factorization (\HQR) yields columns of $Q $ that are orthogonal to a high degree of precision, making these algorithms the weapon of choice in many situations.   Pivoting can be added to \HQR\ to yield \HQR\ with column pivoting (\HQRP).  This topic is covered by standard texts on numerical linear algebra~\cite{GVL3}.

To achieve high performance for dense linear algebra algorithms,
so-called blocked algorithms are employed that cast most computation in terms of matrix-matrix operations supported by the widely used level-3 Basic Linear Algebra Subprograms  (BLAS)~\cite{BLAS3,DDSV} because
such operations can be implemented to achieve very high performance on modern processors via a combination of careful reuse of data in the caches and low level implementation in terms of assembly code or intrinsic vector operations.  Widely used current implementations of the level-3 BLAS are based on
techniques exposed by Goto~\cite{Goto,Goto2} and available in open source libraries including the OpenBLAS~\cite{xianyi2012model} (a fork of the GotoBLAS) and BLIS~\cite{BLIS1}, as well as vendor implementions including AMD's ACML~\cite{ACML}, Intel's MKL~\cite{MKL}, and IBM's ESSL~\cite{ESSL} libraries.

The fundamental problem with the classical approach to \HQRP\ is that only half of the computation can be cast in terms of \gemm, as described in the paper~\cite{QRP:SIAM} that underlies LAPACK's {\tt geqp3} routine~\cite{LAPACK3}.
This means that blocking can only improve performance by, at best, a factor two,  which is inherent from the fact that it must be known how remaining columns will be updated in order to compute  the 2-norms of remaining columns.
Bischof and Quintana-Ort\'{\i} describe in a pair of papers~\cite{RRQRblk1,RRQRblk2} an
attempt to overcome this problem by using so called ``window pivoting'' in combination with \HQR.
While much faster than {\tt geqp3}, this approach is more complicated than
the method proposed in this paper and never made it into LAPACK.

\input Graphs/speedup_sidebyside

The present paper proposes to solve the problem by means of randomized projections.
To describe the idea, suppose that we seek to determine a set of $b$ good pivot columns
in an $m\times n$ matrix $A$. We then draw a Gaussian random matrix $\mtx{G}$ of size
$b\times m$ and form a $b\times n$ \textit{sampling matrix} $\mtx{Y} = \mtx{G}\mtx{A}$. Once $\mtx{Y}$
is available, we execute QRP on this matrix to find the $b$ pivot columns. This computation
is efficient since $\mtx{Y}$ is small compared to $\mtx{A}$ (it has only $b$ rows), and results 
in good pivot choices since the random projection produces a matrix $\mtx{Y}$ that has approximately 
the same linear dependencies between its columns as does $\mtx{A}$.\footnote{To be precise, for linear 
dependencies to be preserved reliably, one needs to perform a very slight amount of over-sampling.
See Section \ref{sec:basicidea} for details.}
With this observation, it becomes easy to block the Householder QR factorization with column
pivoting. At each iteration of the blocked algorithm, we use the randomized sampling approach
to identify a set of $b$ columns that are then moved to the front of the actual matrix,
at which point a regular step of \HQR\ can be used to move the computation forward, optionally
with additional column pivoting only within a narrow panel of the matrix.  Importantly, the
sampling matrix can be cheaply downdated rather than recomputed at each step, allowing the
performance of the proposed algorithm to asymptotically approach that of a standard blocked
\HQR\ implementation that does not pivot columns.
Fig.~\ref{fig:speedup} illustrates the dramatic performance improvements that are realized.
\NoShow{Performance experiments demonstrate the
effectiveness of the proposed approach and the quality of the resulting QR factorization is
also illustrated.}
%A random matrix is used to create a much smaller, representative matrix that we call the sampling marix.
% By computing the QRP of that matrix, %in each iteration of the blocked algorithm columns are identified
% that are swapped to the front of the actual matrix, at which point a regular %step of \HQR\ can be used
% to move the computation forward, optionally with additional column pivoting only within a narrow panel
% of the matrix.  %Importantly, the sampling matrix can be cheaply downdated rather than recomputed at
% each step, allowing the performance of the proposed %algorithm to asymptotically approach that of a
% standard blocked  \HQR\ implementation that down not pivot columns.  Performance experiments demonstrate
% the effectiveness of the proposed approach and the quality of the resulting QR factorization is also illustrated.

%in a fashion
%that was originally described in \cite{2006_martinsson_random1_orig} and later in
%\cite{2007_martinsson_PNAS,2011_martinsson_randomsurvey,2011_martinsson_random1}.

The idea to use randomized sampling to pick blocks of pivot vectors was first published by Martinsson
on ArXiv in May 2015~\cite{2015_blockQR}. A very similar technique was published by Duersch and Gu in
September 2015~\cite{2015_blockQR_ming}, also on ArXiv. The observation that downdating of the sampling
matrix enables the randomized scheme to attain the same asymptotic flop count as classical
\HQRP~was discovered independently by the two groups and was first published in \cite{2015_blockQR_ming}.
More broadly, the idea that one can select a subset of columns of a matrix that forms a good approximate basis for the
column space of the matrix by performing \QRP~on a small matrix whose rows are random linear combinations
of the rows of the original matrix was first described in \cite[Sec.~4.1]{2006_martinsson_random1_orig}
and later elaborated in \cite{2007_martinsson_PNAS,2011_martinsson_randomsurvey,2011_martinsson_random1}.
This problem is closely related to the problem of finding a set of columns of maximal spanning volume
\cite{gu1996}, and to the problem of finding so called \textit{CUR} and \textit{interpolative}
decompositions \cite{2014_martinsson_CUR}.
These ideas tie in to a larger literature on randomized techniques for computing low-rank approximations
of matrices that includes \cite{2004_kannan_vempala,2006_drineas_kannan_mahoney,mahoney2011randomized,2015arXiv150307157M}.
%In particular, \cite{2007_martinsson_PNAS} demonstrates that the randomized method is an efficient and accurate tool for
%finding a set of columns of a matrix that form a good basis for the column space.

%The basic idea of using randomized projections to
%efficiently identify a set of good pivot columns can be traced to \cite{2007_martinsson_PNAS}. In
%\cite{2007_martinsson_PNAS}, the idea is used to find a set of columns of a matrix with close to
%maximal spanning volume, which is a problem that is very closely related to the problem of finding
%good pivot columns \cite{gu1996}. The idea of \cite{2007_martinsson_PNAS} was subsequently refined in
%\cite[Sec.~5.2]{2011_martinsson_randomsurvey}, and used in \cite{2014_martinsson_CUR} to the problem
%of computing so called \textit{CUR} and \textit{interpolative} decompositions, which are closely
%related to the problem of finding rank-revealing QR factorizations.

This paper describes a practical implementation of the proposed method that can be incorporated
in libraries like LAPACK and {\tt libflame}~\cite{libflame_ref,CiSE09}.
Implementation details that are important for attaining high practical performance are
described to enable readers to reproduce and extend the ideas.
The paper provides a cost analysis that shows that asymptotically the number of floating point
operations approaches that of \HQR\ without pivoting while most computation is cast in terms of
matrix-matrix multiplication like the corresponding blocked \HQR\  without pivoting.
It reports unprecedented performance for pivoted QR factorization on current architectures
and provides empirical quality results.
Importantly, the implementation is made available for use by the computational science community
under an Open Source license.
The conclusion discusses how these results pave the way for future opportunities.

The paper is organized as follows: Section 2 lists some standard facts about Householder
reflectors and pivoted QR factorizations that we need in the presentation. Section 3 describes
how the classical Householder QR factorization algorithm can be blocked by using randomization in the pivot
selection step. Section 4 reports the results from numerical experiments investigating the speed
of the algorithm and the quality of the pivoting selection strategy. Sections 5 summarizes the
key results and discusses future work. Section 6 describes publicly available software that
implements the techniques presented. 

%% file: Graphs/speedup_sidebyside.tex
%\begin{wrapfigure}{r}{1.0\textwidth}
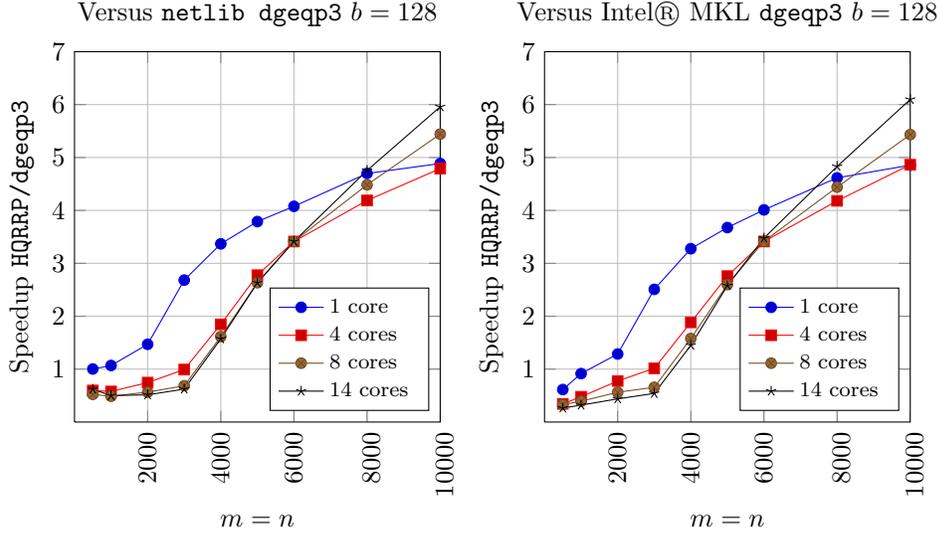
\begin{figure}[tb]
\begin{center}
\begin{tikzpicture}
  \begin{axis}[
    title= {Versus {\tt netlib dgeqp3} $ b=128$},
               xtick={%
%                1,
                 2,
%                3,
                 4,
%                5,
                 6,
%                7,
                 8,
%                9,
                10},
               xticklabels={%
%                 \mbox{\begin{sideways}1000\end{sideways}},
                 \mbox{\begin{sideways}2000\end{sideways}},
%                 \mbox{\begin{sideways}3000\end{sideways}},
                 \mbox{\begin{sideways}4000\end{sideways}},
%                 \mbox{\begin{sideways}5000\end{sideways}},
                 \mbox{\begin{sideways}6000\end{sideways}},
%                 \mbox{ \begin{sideways}7000\end{sideways}},
                 \mbox{ \begin{sideways}8000\end{sideways}},
%                 \mbox{ \begin{sideways}9000\end{sideways}},
                 \mbox{ \begin{sideways}10000\end{sideways}}
               },
               ytick={1,2,3,4,5,6,7}
               ,
               xmajorgrids,ymajorgrids,
               width=0.495\textwidth,
               height=0.5\textwidth,
               xlabel=\mbox{$m =n$},
               ylabel={Speedup \tt HQRRP/dgeqp3},
               xmin=0,xmax=10,ymin=0,ymax=7,
               legend style={legend pos=south east},
%               legend pos=outer north east,
               clip=true],
%,cycle list name=mark list],
],
% qrrp with downdate
  \addplot coordinates {
(   500/1000,   4.07 /  4.06 )
(  1000/1000,   7.35 /  6.88 )
(  2000/1000,  12.15 /  8.27 )
(  3000/1000,  15.13 /  5.64 )
(  4000/1000,  17.11 /  5.08 )
(  5000/1000,  18.83 /  4.97 )
(  6000/1000,  20.10 /  4.93 )
(  8000/1000,  21.99 /  4.68 )
( 10000/1000,  23.21 /  4.75 )
  };
  \addplot coordinates {
(   500/1000,   3.33 /  5.59 )
(  1000/1000,   7.60 / 13.16 )
(  2000/1000,  15.94 / 21.38 )
(  3000/1000,  23.98 / 24.17 )
(  4000/1000,  30.74 / 16.65 )
(  5000/1000,  37.42 / 13.48 )
(  6000/1000,  43.13 / 12.64 )
(  8000/1000,  52.56 / 12.55 )
( 10000/1000,  59.69 / 12.46 )
  };
  \addplot coordinates {
(   500 / 1000,   3.56 /  6.81 )
(  1000 / 1000,   8.06 / 16.55 )
(  2000 / 1000,  18.07 / 31.98 )
(  3000 / 1000,  25.76 / 37.58 )
(  4000 / 1000,  34.31 / 21.24 )
(  5000 / 1000,  43.35 / 16.46 )
(  6000 / 1000,  51.43 / 15.07 )
(  8000 / 1000,  66.50 / 14.83 )
( 10000 / 1000,  78.86 / 14.49 )
  };
  \addplot coordinates {
(   500 / 1000,   2.75 /  4.43 )
(  1000 / 1000,   6.95 / 13.91 )
(  2000 / 1000,  17.66 / 34.35 )
(  3000 / 1000,  25.52 / 41.01 )
(  4000 / 1000,  35.32 / 22.46 )
(  5000 / 1000,  45.18 / 17.18 )
(  6000 / 1000,  54.75 / 16.04 )
(  8000 / 1000,  73.81 / 15.52 )
( 10000 / 1000,  90.87 / 15.25 )
  };
  \legend{%
    \begin{minipage}{0.45in}\footnotesize  1 core\end{minipage},
    \begin{minipage}{0.45in}\footnotesize  4 cores\end{minipage},
    \begin{minipage}{0.45in}\footnotesize  8 cores\end{minipage},
    \begin{minipage}{0.45in}\footnotesize 14 cores\end{minipage},
  }
  \end{axis}
\end{tikzpicture}
\begin{tikzpicture}
  \begin{axis}[
    title= {Versus Intel\textregistered\ MKL {\tt dgeqp3} $ b=128$},
               xtick={%
%                1,
                 2,
%                3,
                 4,
%                5,
                 6,
%                7,
                 8,
%                9,
                10},
               xticklabels={%
%                 \mbox{\begin{sideways}1000\end{sideways}},
                 \mbox{\begin{sideways}2000\end{sideways}},
%                 \mbox{\begin{sideways}3000\end{sideways}},
                 \mbox{\begin{sideways}4000\end{sideways}},
%                 \mbox{\begin{sideways}5000\end{sideways}},
                 \mbox{\begin{sideways}6000\end{sideways}},
%                 \mbox{ \begin{sideways}7000\end{sideways}},
                 \mbox{ \begin{sideways}8000\end{sideways}},
%                 \mbox{ \begin{sideways}9000\end{sideways}},
                 \mbox{ \begin{sideways}10000\end{sideways}}
               },
               ytick={1,2,3,4,5,6,7}
               ,
               xmajorgrids,ymajorgrids,
               width=0.495\textwidth,
               height=0.5\textwidth,
               xlabel=\mbox{$m =n$},
               ylabel={Speedup \tt HQRRP/dgeqp3},
               xmin=0,xmax=10,ymin=0,ymax=7,
               legend style={legend pos=south east},
%               legend pos=outer north east,
               clip=true],
%,cycle list name=mark list],
],
% qrrp with downdate
  \addplot coordinates {
( 500/1000,     4.07 /  6.63 )
( 1000/1000,    7.38 /  8.06 )
( 2000/1000,   12.04 /  9.36 )
( 3000/1000,   14.97 /  5.97 )
( 4000/1000,   17.10 /  5.22 )
( 5000/1000,   18.71 /  5.09 )
( 6000/1000,   20.01 /  4.99 )
( 8000/1000,   21.83 /  4.73 )
( 10000/1000,  23.21 /  4.78 )
  };
    \addplot coordinates {
(   500/1000,   3.29 /  9.59 )
(  1000/1000,   7.51 / 15.72 )
(  2000/1000,  15.94 / 20.56 )
(  3000/1000,  23.90 / 23.56 )
(  4000/1000,  30.74 / 16.32 )
(  5000/1000,  37.35 / 13.54 )
(  6000/1000,  42.82 / 12.55 )
(  8000/1000,  52.09 / 12.46 )
( 10000/1000,  60.02 / 12.34 )
  };
  \addplot coordinates {
(   500 / 1000,   3.53 / 10.71 )
(  1000 / 1000,   8.13 / 20.68 )
(  2000 / 1000,  18.07 / 32.10 )
(  3000 / 1000,  25.80 / 39.30 )
(  4000 / 1000,  34.43 / 21.81 )
(  5000 / 1000,  43.28 / 16.66 )
(  6000 / 1000,  51.58 / 15.10 )
(  8000 / 1000,  66.33 / 14.94 )
( 10000 / 1000,  79.00 / 14.54 )
  };
  \addplot coordinates {
(   500 / 1000,   2.74 / 10.37 )
(  1000 / 1000,   7.01 / 21.95 )
(  2000 / 1000,  17.45 / 39.83 )
(  3000 / 1000,  25.45 / 47.08 )
(  4000 / 1000,  35.14 / 24.23 )
(  5000 / 1000,  45.19 / 17.53 )
(  6000 / 1000,  55.28 / 15.90 )
(  8000 / 1000,  73.43 / 15.20 )
( 10000 / 1000,  90.61 / 14.86 )
  };
  \legend{%
    \begin{minipage}{0.45in}\footnotesize  1 core\end{minipage},
    \begin{minipage}{0.45in}\footnotesize  4 cores\end{minipage},
    \begin{minipage}{0.45in}\footnotesize  8 cores\end{minipage},
    \begin{minipage}{0.45in}\footnotesize 14 cores\end{minipage},
  }
  \end{axis}
\end{tikzpicture}
\end{center}

\caption{Speedup of new blocked Householder QR factorization with
  randomized column pivoting ({\tt HQRRP}) relative to LAPACK's
  faster routine ({\tt dgeqp3}) on a 14-core Intel Xeon E5-2695 v3, see Section~\ref{sec:performance} for details.
}
\label{fig:speedup}

%\end{wrapfigure}
\end{figure}

%% file: 02HQR.tex
In this section, we briefly review the state-of-the art regarding
Householder factorization based on Householder transformations (\HQR).
Throughout, we use the FLAME notation for representing dense linear algebra algorithms~\cite{inverse-siam,FLAME}.
In particular, for any matrix $\mtx{X}$, we let $m(\mtx{X})$ and $n(\mtx{X})$ denote the number of
rows and columns of $\mtx{X}$, respectively.

\subsection{Householder transformations (reflectors)}

A standard topic in numerical linear algebra is the concept of a
reflector, also known as a Householder transformation~\cite{GVL3}.
The review in this subsection
follows~\cite{Joffrain:2006:AHT:1141885.1141886} in which a similar notation is also employed.

Given a nonzero vector $ u \in \Cn $, the matrix $ H( u ) = I - \frac{1}{\tau} u u^H $ with $ \tau = \frac{u^H u}{2} $ has the property that it reflects a vector to which it is applied with respect to the subspace orthogonal to $ u $.
Given a vector $ x $, the vector $ u $ and scalar $ \tau $ can be chosen so that $ H(u) x $ equals a multiple of $ e_0 $,
the first column of the identiy matrix.   Furthermore, $ u$ can be normalized so that its first element equals one.

In our discussions,
given a vector $ x = \FlaTwoByOneSingleLine{ \chi_1 }{ x_2 } $,
the function
$
\left[
\FlaTwoByOneSingleLine
  {\rho}
  { u_2},
\tau \right] \becomes
\housev \left(
\FlaTwoByOneSingleLine
{ \chi_1 }
{ x_2 }  \right)
$
computes
the vector
$ u = \FlaTwoByOneSingleLine{1}{u_2} $ and $ \tau = \frac{u^H u}{2} $
so that $ H( u ) x = \rho e_0 $,

\subsection{Unblocked Householder QR factorization}

\input HQRFormTunb

A standard unblocked algorithm for \HQR~of a  given matrix $ A \in \Cmxn $, typeset using the FLAME notation, is given
in Fig.~\ref{fig:HQR_unb_blk}~(left).
The body of the loop computes
\[
\left[
\FlaTwoByOneSingleLine
  {\rho_{11}}
  { u_{21}},
\tau_{11} \right] \becomes
\housev \left(
\FlaTwoByOneSingleLine
{ \alpha_{11} }
{ a_{21} }  \right),
\]
which overwrites $ a_{11} $ with $ \rho_{11} $ and $ a_{21} $ with $ u_{21}
$, after which the remainder of $ A $ is updated by
\[
\left(
\begin{array}{c}
a_{12}^T \\ \hline
A_{22}
\end{array}
\right) \becomes
\left(
I -
\frac{1}{\tau_{11}}
\left(
\begin{array}{c}
1 \\ \hline
u_{21}
\end{array} \right)
\left(
\begin{array}{c}
1 \\ \hline
u_{21}
\end{array}
\right)^H
\right)
\left(
\begin{array}{c}
a_{12}^T \\ \hline
A_{22}
\end{array}
\right).
\]
%%% via the steps
%%% \begin{itemize}
%%% \item
%%% $ w_{12}^T \becomes ( a_{12}^T + a_{21}^H A_{22} )/\tau_1 $
%%% \item
%%% $ \FlaTwoByOneSingleLine
%%% { a_{12}^T }
%%% { A_{22} }
%%% \becomes
%%% \FlaTwoByOneSingleLine
%%% { a_{12}^T - { w_{12}^T } }
%%% { A_{22} - { a_{21} } { w_{12}^T } }
%%% $
%%% \end{itemize}
Upon completion, the (Householder) vectors that define the Householder transformations have overwritten the elements in which they introduced zeroes, and the upper triangular part of $ A $ contains $ R $.  How the matrix $ T $ fits into the picture will become clear next.

\subsection{The UT transform: Accumulating Householder transformations}

Given $ A \in \mathbb{C}^{n \times b} $, let $ U $ contain the
Householder vectors computed during the \HQR\ of $ A $.
Let us assume that $ H(u_{b-1}) \cdots H(u_1) H(u_0) A = R $.  Then there exists an upper triangular matrix so that
$ I - U T^{-H} U^H = H(u_{b-1}) \cdots H(u_1) H(u_0) $.
The desired matrix $ T $ equals the strictly upper triangular part of $ U^H U $ with the diagonal elements equal to $ \tau_0, \ldots  , \tau_{b-1} $.  The matrix $ T $ can be computed during the unblocked \HQR, as indicated in Fig.~\ref{fig:HQR_unb_blk}~(left).
In~\cite{Joffrain:2006:AHT:1141885.1141886}, the transformation $ I - U T^{-1} U^H $ that equals
the accumulated Householder transformations is called the {\em UT transform}. \pgm{The UT transform is conceptually
related to the more familiar WY transform~\cite{BiVL87} and compact WY transform~\cite{ScVL89}, see
\cite{Joffrain:2006:AHT:1141885.1141886} for details on how the different representations relate to one another.}

\subsection{A blocked QR Householder factorization algorithm}

\input HQR_unb_blk

A blocked algorithm for HQR that exploits the insights that resulted in the UT
transform can now be described as follows.
Partition
\[
A \rightarrow
\FlaTwoByTwoSingleLine
  { A_{11} }{ A_{12} }
  { A_{21} }{ A_{22} }
\]
where $ A_{11} $ is $ b \times b $.
We can use the unblocked algorithm in Fig.~\ref{fig:HQR_unb_blk}~(left) to
factor the panel
$
\FlaTwoByOneSingleLine
  { A_{11} }
  { A_{21} }
$, creating matrix $ T_{11} $ as a side effect.
Now we need to also apply the UT transform to the rest
of the columns:
\begin{eqnarray*}
\FlaTwoByOneSingleLine
  { A_{12} }
  { A_{22} }
&:=&
\left(
I -
\FlaTwoByOneSingleLine
  { U_{11} }
  { U_{21} }
T_{11}^{-1}
\FlaTwoByOneSingleLine
  { U_{11} }
  { U_{21} }^H
\right)^H
\FlaTwoByOneSingleLine
  { A_{12} }
  { A_{22} }
%%% T=&
%%% \FlaTwoByOneSingleLine
%%%   { A_{12} }
%%%  T{ A_{22} }
%%% T
%%% \FlaTwoByOneSingleLine
%%%   { U_{11} }
%%%   { U_{21} }
%%%   W_{12} \\
=
\FlaTwoByOneSingleLine
  { A_{12} - U_{11} W_{12} }
  { A_{22} - U_{21} W_{12} }
,
\end{eqnarray*}
where
$
W_{12} =
T_{11}^{-H} ( U_{11}^H A_{12} + U_{21}^H A_{22} )
$.
This motivates the blocked HQR algorithm in Fig.~\ref{fig:HQR_unb_blk}~(right) which we will refer to as {\sc HQR\_blk}.

The benefit of the blocked algorithm is that it casts most computation in terms of
the computations $ U_{21}^H A_{22} $ (row panel times matrix multiply) and
$ A_{22} - U_{21} W_{12} $ (rank-$b$ update).  Such matrix-matrix multiplications can attain high performance by amortizing data movement between memory layers.

These insights form the basis for the LAPACK routine {\sc geqrf}
(except that it uses a compact WY transform instead of the UT transform).

\subsection{Householder QR factorization with column pivoting}

\NoShow{\input HQRP_unb_blk}

An unblocked (rank-revealing) Householder QR factorization with column
pivoting ({\sc HQRP}) swaps the column of $ A_{BR} $ with largest
2-norm with the first column of that matrix at the top of the loop
body.  As a result, the diagonal elements of matrix $ R $ are ordered
from largest to smallest in magnitude, which, for example, allows the
resulting QR factorization to be used to identify a high quality
approximate low-rank orthonormal basis for the column space of $ A $.
(To be precise, column pivoted QR returns a high quality basis in most
cases but may produce strongly sub-optimal results in rare situations.
For details, see \cite{1966_kahan_NLA,gu1996}, and the description of
``Matrix 4'' in Section \ref{sec:pivotquality}.)

The fundamental problem with the best known algorithm for HQRP, which
underlies LAPACK's routine {\tt geqp3},
is that it only casts half of the computation in terms of
matrix-matrix multiplication~\cite{QRP:SIAM}.  The unblocked algorithm
called from the blocked algorithm operates on the entire ``remaining
matrix'' ($ A_{BR} $ in the blocked algorithm), computes $ b $ more
Householder transforms and $ b $ more rows of $ R $, computes the
matrix $ W_2 $, and returns the information about how columns were
swapped.  In the blocked algorithm itself, only the update $ A_{22} -
A_{21} W_2 $ remains to be performed.  When only half the computation
can be cast in terms of matrix-matrix multiplication, the resulting
blocked algorithm is only about twice as fast as the unblocked
algorithm.

%% file: HQRFormTunb.tex
\resetsteps      % Reset all the commands to create a blank worksheet  

% Define the operation to be computed

\renewcommand{\routinename}{%
\mbox{~~~~~~~~~$ \left[ A, T \right] := \mbox{\sc HQR\_unb}( A, T ) $}
}

% Step 3: Loop-guard 

\renewcommand{\guard}{
  m( A_{TL} ) < m( A )
}

% Step 4: Define Initialize 

\renewcommand{\partitionings}{
  $
  A \rightarrow
  \FlaTwoByTwo{A_{TL}}{A_{TR}}
              {A_{BL}}{A_{BR}}
  $
,
\\
~~~  $
  T \rightarrow
  \FlaTwoByTwo{T_{TL}}{T_{TR}}
              {0}{T_{BR}}
  $
}

\renewcommand{\partitionsizes}{
$ A_{TL} $ is $ 0 \times 0 $,
$ T_{TL} $ is $ 0 \times 0 $
}

% Step 5a: Repartition the operands 

\renewcommand{\repartitionings}{
$  \FlaTwoByTwo{A_{TL}}{A_{TR}}
              {A_{BL}}{A_{BR}}
  \rightarrow
  \FlaThreeByThreeBR{A_{00}}{a_{01}}{A_{02}}
                    {a_{10}^T}{\alpha_{11}}{a_{12}^T}
                    {A_{20}}{a_{21}}{A_{22}}
$,
\\
$  \FlaTwoByTwo{T_{TL}}{T_{TR}}
              {0}{T_{BR}}
  \rightarrow
  \FlaThreeByThreeBR{T_{00}}{t_{01}}{T_{02}}
                    {0}{\tau_{11}}{t_{12}^T}
                    {0}{0}{T_{22}}
$}

\renewcommand{\repartitionsizes}{
  $ \alpha_{11} $ is $ 1 \times 1 $,
  $ \tau_{11} $ is $ 1 \times 1 $}

% Step 5b: Move the double lines 

\renewcommand{\moveboundaries}{
$  \FlaTwoByTwo{A_{TL}}{A_{TR}}
              {A_{BL}}{A_{BR}}
  \leftarrow
  \FlaThreeByThreeTL{A_{00}}{a_{01}}{A_{02}}
                    {a_{10}^T}{\alpha_{11}}{a_{12}^T}
                    {A_{20}}{a_{21}}{A_{22}}
$,
\\
$  \FlaTwoByTwo{T_{TL}}{T_{TR}}
              {0}{T_{BR}}
  \leftarrow
  \FlaThreeByThreeTL{T_{00}}{t_{01}}{T_{02}}
                    {0}{\tau_{11}}{t_{12}^T}
                    {0}{0}{T_{22}}
$}

% Step 8: Insert the updates required to change the 
%         state from that given in Step 6 to that given in Step 7
% Note: The below needs editing!!!

\renewcommand{\update}{
$
\begin{array}{l}
\left[
\FlaTwoByOneSingleLine
{ \alpha_{11} }
{ a_{21}}, \tau_{11}
\right]
\becomes
%%% \left[
%%% \FlaTwoByOneSingleLine
%%% { \rho_{11} }
%%% { u_{21}}, \tau_{11}
%%% \right]
%%% =
\housev
\left( \begin{array}{c}
\alpha_{11} \\ \hline
a_{21}
\end{array}
\right) \\
w_{12}^T \becomes ( a_{12}^T + a_{21}^H A_{22} )/\tau_{11} 
\\
\FlaTwoByOneSingleLine
{ a_{12}^T }
{ A_{22} }
\becomes 
\FlaTwoByOneSingleLine
{ a_{12}^T - { w_{12}^T } }
{ A_{22} - { a_{21} } { w_{12}^T } }
\\
t_{01} := ( a_{10}^T )^H + A_{20}^H a_{21} \\[0.33in]
\end{array}
$ 
}

%%% \begin{figure}[tb!]
%%% \FlaAlgorithm
%%% \caption{Unblocked Householder transformation based QR factorization
%%%   merged with the computation of $ T $ for the UT transform.}
%%% \label{fig:HQRunb}
%%% \end{figure}

%% file: HQR_unb_blk.tex
\begin{figure}[tb!]

\begin{center}
\begin{tabular}{@{}c c}
\begin{minipage}{0.47\textwidth}
\setlength{\arraycolsep}{2pt}

\input HQRFormTunb
\FlaAlgorithmNarrow
\end{minipage}
&
\begin{minipage}{0.45\textwidth}
\setlength{\arraycolsep}{2pt}

\input HQRblk_var2
\FlaAlgorithmNarrow
\end{minipage}
\end{tabular}
\end{center}
\caption{
Left: Unblocked Householder transformation based QR factorization
merged with the computation of $ T $ for the UT transform.
Right:
Blocked Householder transformation based QR factorization.  In this algorithm, $ U_{11} $ is the unit lower triangular matrix stored below the diagonal of $ A_{11} $ and $ U_{21} $ is stored in $ A_{21} $.}
\label{fig:HQR_unb_blk}
\end{figure}

%% file: HQRblk_var2.tex
\resetsteps      % Reset all the commands to create a blank worksheet  

% Define the operation to be computed

\renewcommand{\routinename}{%
\mbox{~~~~~~~~~$  \left[ A, T \right] := \mbox{\sc HQR\_blk}( A, T ) $
} 
}

% Step 3: Loop-guard 

\renewcommand{\guard}{
  m( A_{TL} ) < m( A )
}

% Step 4: Define Initialize 

\renewcommand{\partitionings}{
  $
  A \rightarrow
  \FlaTwoByTwo{A_{TL}}{A_{TR}}
              {A_{BL}}{A_{BR}}
  $
,\\
  $
  T \rightarrow
  \FlaTwoByOne{T_{T}}
              {T_{B}}
  $
}

\renewcommand{\partitionsizes}{
$ A_{TL} $ is $ 0 \times 0 $,
$ T_{T} $ has $ 0 $ rows
}

% Step 5a: Repartition the operands 

\renewcommand{\blocksize}{b}

\renewcommand{\repartitionings}{
$  \FlaTwoByTwo{A_{TL}}{A_{TR}}
              {A_{BL}}{A_{BR}}
  \rightarrow
  \FlaThreeByThreeBR{A_{00}}{A_{01}}{A_{02}}
                    {A_{10}}{A_{11}}{A_{12}}
                    {A_{20}}{A_{21}}{A_{22}}
$,
\\
$  \FlaTwoByOne{ T_T }
              { T_B }
\rightarrow
  \FlaThreeByOneB{T_0}
                 {T_1}
                 {T_2}
$
}

\renewcommand{\repartitionsizes}{
  $ A_{11} $ is $ b \times b $,
  $ T_1 $ has $ b $ rows}

% Step 5b: Move the double lines 

\renewcommand{\moveboundaries}{
$  \FlaTwoByTwo{A_{TL}}{A_{TR}}
              {A_{BL}}{A_{BR}}
  \leftarrow
  \FlaThreeByThreeTL{A_{00}}{A_{01}}{A_{02}}
                    {A_{10}}{A_{11}}{A_{12}}
                    {A_{20}}{A_{21}}{A_{22}}
$,
\\
$  \FlaTwoByOne{ T_T }
               { T_B }
\leftarrow
  \FlaThreeByOneT{T_0}
                 {T_1}
                 {T_2}
$
}

% Step 8: Insert the updates required to change the 
%         state from that given in Step 6 to that given in Step 7
% Note: The below needs editing!!!

\renewcommand{\update}{
$
  \begin{array}{l}
[
\FlaTwoByOneSingleLine
  { A_{11} }
  { A_{21} }, T_{1} ] 
:= 
\\
~~~~~ \mbox{\sc HQR\_unb}(
\FlaTwoByOneSingleLine
  { A_{11} }
  { A_{21} }, T_1 ) \\
W_{12} := 
T_{1}^{-H} ( U_{11}^H A_{12} + U_{21}^H A_{22} ) \\
\FlaTwoByOneSingleLine
  { A_{12} }
  { A_{22} }
:= 
\FlaTwoByOneSingleLine
  { A_{12} - U_{11} W_{12} }
  { A_{22} - U_{21} W_{12} }
\\[-0.12in]
~
\end{array}
$
}

%%% \begin{figure}[tb!]
%%% \FlaAlgorithm
%%% \caption{Blocked Householder transformation based QR
%%% factorization.  In this algorithm, $ U_{11} $ is the unit lower
%%% triangular matrix stored below the diagonal of $ A_{11} $ and
%%% $ U_{21} $ is stored in $ A_{21} $.}
%%% \label{fig:HQRblk}
%%% \end{figure}

%% file: HQRP_unb_blk.tex
\begin{figure}[tb!]

\begin{center}
\begin{tabular}{@{}c c}
\begin{minipage}{0.47\textwidth}
\input RRQR_unb_var3
\FlaAlgorithmNarrow
\end{minipage}
&
\begin{minipage}{0.47\textwidth}
\input RRQR_blk
\FlaAlgorithmWithInitNarrow
\end{minipage}
\end{tabular}
\end{center}
\caption{
{\bf Left:} Unblocked HQRP. 
{\bf Right:} Blocked HQRP.  
}
\label{fig:HQRP_unb_blk}
\end{figure}

%% file: RRQR_unb_var3.tex
\resetsteps

\renewcommand{\initialize}{
\phantom{ v := \mbox{ComputeWeights}}
}

\renewcommand{\routinename}{%
\hspace{-0.75in}\begin{minipage}[t]{2.75in}
~\\
$\left[ A, T, s, v, W \right] = \mbox{\sc 
    HQRP\_unb
}( A, T, s, v, W, r ) $
\end{minipage}}

\renewcommand{\guard}{n( A_{TL} ) < r}

\renewcommand{\partitionings}{
\footnotesize
$ 
A \rightarrow 
\FlaTwoByTwo{ A_{TL} } {A_{TR}}
              { A_{BL} } { A_{BR}   }
$
, \\
~~~ $
T \rightarrow 
\FlaTwoByTwo{ T_{TL} } {T_{TR}}
              { T_{BL} } { T_{BR}   }
$
,
$
s \rightarrow 
\FlaTwoByOne{ s_{T} }
            { s_{B} }
$
, \\
~~~ $
v \rightarrow 
\FlaTwoByOne{ v_{T} }
            { v_{B} }
$,
$ 
W \rightarrow 
\FlaTwoByTwo{ W_{TL} } {W_{TR}}
              { W_{BL} } { W_{BR}   }
$
}

\renewcommand{\partitionsizes}{
$ A_{TL} $, $ T_{TL} $, $ W_{TL} $ are $ 0 \times 0 $, \\
$ v_T $ has $ 0 $ elements}

\renewcommand{\repartitionings}{
\footnotesize 
$
\FlaTwoByTwo{ A_{TL} } {A_{TR}}
              { A_{BL} } { A_{BR}   } \rightarrow 
\FlaThreeByThreeBR{ A_{00} }{ a_{01} }{ A_{02} }
              { a_{10}^T }{ \alpha_{11} }{ a_{12}^T }
              { A_{20} }{ a_{21} }{ A_{22} }
$
, \\
$
\FlaTwoByTwo{ T_{TL} } {T_{TR}}
              { T_{BL} } { T_{BR}   }
\rightarrow 
\FlaThreeByThreeBR{ T_{00} }{ t_{01} }{ T_{02} }
              { 0 }{ \tau_{11} }{ t_{12}^T }
              { 0 }{ 0 }{ T_{22} }
$
,\\
$
\FlaTwoByOne{ s_{T} }
            { s_{B} } \rightarrow 
\FlaThreeByOneB{ s_{0} }
              { \sigma_{1} }
              { s_{2} }
$
, 
$
\FlaTwoByOne{ v_{T} }
            { v_{B} } \rightarrow 
\FlaThreeByOneB{ v_{0} }
              { \nu_{1} }
              { v_{2} }
$,\\
$
\FlaTwoByTwo{ W_{TL} } {W_{TR}}
              { W_{BL} } { W_{BR}   } \rightarrow 
\FlaThreeByThreeBR{ W_{00} }{ w_{01} }{ W_{02} }
              { w_{10}^T }{ \omega_{11} }{ w_{12}^T }
              { W_{20} }{ w_{21} }{ W_{22} }
$
}
\renewcommand{\repartitionsizes}{
$ \alpha_{11} $, $ \tau_{11} $, $ \cdots $ are  scalars
}

\renewcommand{\moveboundaries}{
\NoShow{
$
\FlaTwoByTwo{ A_{TL} } {A_{TR}}
              { A_{BL} } { A_{BR}   } \leftarrow
\FlaThreeByThreeTL{ A_{00} }{ a_{01} }{ A_{02} }
              { a_{10}^T }{ \alpha_{11} }{ a_{12}^T }
              { A_{20} }{ a_{21} }{ A_{22} }
$
,
$
\FlaTwoByTwo{ T_{TL} } {T_{TR}}
              { T_{BL} } { T_{BR}   }
\leftarrow
\FlaThreeByThreeTL{ T_{00} }{ t_{01} }{ T_{02} }
              { 0 }{ \tau_{11} }{ t_{12}^T }
              { 0 }{ 0 }{ T_{22} }
$
}
$ \cdots $
}

\renewcommand{\update}{
$
\begin{array}{l}
[
\left( \begin{array}{c}
\nu_1 \\ \hline 
v_2 
\end{array}
\right),
\sigma_1 ] = \mbox{\sc DeterminePivot}( 
\left( \begin{array}{c}
\nu_1 \\ \hline 
v_2 
\end{array}
\right) ) \\
\left( \begin{array}{c | c }
a_{01} & A_{02} \\ \whline 
\alpha_{11} & a_{12}^T \\ \hline 
a_{21} & A_{22} 
\end{array}
\right) 
\becomes
\mbox{\sc Swap}(
\sigma_1,
\left( \begin{array}{c | c}
a_{01} & A_{02} \\ \whline 
\alpha_{11} & a_{12}^T \\ \hline 
a_{21} & A_{22} 
\end{array}
\right) 
)
\\
\alpha_{11} := \alpha_{11} - a_{10}^T w_{01}
\\
\color{red}  
a_{21} := a_{21} - A_{20} w_{01}
\\
a_{12}^T := a_{12}^T - a_{10}^T W_{02}
\\
\left[
\FlaTwoByOneSingleLine
{ \alpha_{11} }
{ a_{21}}, \tau_{11}
\right]
\becomes
%%% \left[
%%% \FlaTwoByOneSingleLine
%%% { \rho_{11} }
%%% { u_{21}}, \tau_{11}
%%% \right]
%%% =
\housev
\left( \begin{array}{c}
\alpha_{11} \\ \hline
a_{21}
\end{array}
\right) \\
w_{12}^T \becomes ( a_{12}^T + a_{21}^H A_{22}-  ({\color{red}  
a_{21}^H A_{20}}) W_{02}
  ) /\tau_{11} \\
 a_{12}^T := a_{12}^T - w_{12}^T
\\
v_2 = \mbox{\sc DowndateWeights}( v_2, a_{12} )
\end{array}
$
}

%%% \begin{figure}[tbp]
%%% %%% \begin{algorithm}
%%% %%% $\left[ A, t \right] \becomes \URt( A )$
%%% %%% \end{algorithm}
%%% \footnotesize
%%% \FlaAlgorithm
%%% \caption{{\color{red} update caption.  Also, $ T $ needs updated.} Unblocked HQRP that up%%% dates $ r $ rows and columns of $ 
%%% $, but leaves the final trailing matrix $ A_{BR} $ prestine.
%%% Here $ v $ has already been initialized before calling
%%% the routine so that it can be called from a blocked algorithm.
%%% }
%%% \label{fig:RRQR_unb_var3}
%%% \end{figure}

%% file: RRQR_blk.tex
\resetsteps      % Reset all the commands to create a blank worksheet  

% Define the operation to be computed

\renewcommand{\routinename}{ 
\hspace{-0.75in}\begin{minipage}[t]{2.75in}
~\\
$\left[ A, T, s \right] := \mbox{\sc
    HQRP\_blk}( A, T, s , r)   $
\end{minipage}}

% Step 3: Loop-guard 

\renewcommand{\guard}{
  n( A_{TL} ) < r )
}

% Step 4: Define Initialize 

\renewcommand{\partitionings}{
  $
  A \rightarrow
  \FlaTwoByTwo{A_{TL}}{A_{TR}}
              {A_{BL}}{A_{BR}}
  $
, \\
~~~
  $
  T \rightarrow
  \FlaTwoByOne{T_{T}}
              {T_{B}}
  $
,
  $
  s \rightarrow
  \FlaTwoByOne{s_{T}}
              {s_{B}}
  $
, \\
~~~
  $
  v \rightarrow
  \FlaTwoByOne{v_{T}}
              {v_{B}}
  $
,
  $
  W \rightarrow
  \FlaOneByTwo{W_L}{W_R}
  $
}

\renewcommand{\partitionsizes}{
$ A_{TL} $ is $ 0 \times 0 $;
$ W_L $ has $ 0 $ cols; \\
$ T_{T} $,
$ s_{T} $, and
$ v_{T} $ have $ 0 $ rows
}

% Step 5a: Repartition the operands 

\renewcommand{\blocksize}{b}

\renewcommand{\repartitionings}{
$  \FlaTwoByTwo{A_{TL}}{A_{TR}}
              {A_{BL}}{A_{BR}}
  \rightarrow
  \FlaThreeByThreeBR{A_{00}}{A_{01}}{A_{02}}
                    {A_{10}}{A_{11}}{A_{12}}
                    {A_{20}}{A_{21}}{A_{22}}
$, \\
~~~
$  \FlaTwoByOne{ T_T }
              { T_B }
\rightarrow
  \FlaThreeByOneB{T_0}
                 {T_1}
                 {T_2}
$
,  \\
$  \FlaTwoByOne{ s_T }
              { s_B }
\rightarrow
  \FlaThreeByOneB{s_0}
                 {s_1}
                 {s_2}
$
, 
$  \FlaTwoByOne{ v_T }
              { v_B }
\rightarrow
  \FlaThreeByOneB{v_0}
                 {v_1}
                 {v_2}
$
, \\
$  \FlaOneByTwo{W_L}{W_R}
\rightarrow  \FlaOneByThreeR{W_0}{W_1}{W_2}
$
}

\renewcommand{\repartitionsizes}{
  $ A_{11} $ is $ b \times b $,
  $ W_1 $ has $b$ cols,
  $ T_1 $,  $ s_1 $,  $ v_1 $ have $ b $ rows
}

% Step 5b: Move the double lines 

\renewcommand{\moveboundaries}{
$ \cdots $
\NoShow{
$  \FlaTwoByTwo{A_{TL}}{A_{TR}}
              {A_{BL}}{A_{BR}}
  \leftarrow
  \FlaThreeByThreeTL{A_{00}}{A_{01}}{A_{02}}
                    {A_{10}}{A_{11}}{A_{12}}
                    {A_{20}}{A_{21}}{A_{22}}
$,
$  \FlaTwoByOne{ T_T }
               { T_B }
\leftarrow
  \FlaThreeByOneT{T_0}
                 {T_1}
                 {T_2}
$
, \\
$  \FlaTwoByOne{ s_T }
               { s_B }
\leftarrow
  \FlaThreeByOneT{s_0}
                 {s_1}
                 {s_2}
$
,
$  \FlaTwoByOne{ v_T }
               { v_B }
\leftarrow
  \FlaThreeByOneT{v_0}
                 {v_1}
                 {v_2}
$
,
$  \FlaOneByTwo{W_L}{W_R}
\leftarrow  \FlaOneByThreeL{W_0}{W_1}{W_2}
$
}
}

% Step 8: Insert the updates required to change the 
%         state from that given in Step 6 to that given in Step 7
% Note: The below needs editing!!!

\renewcommand{\update}{
\normalsize
\setlength{\arraycolsep}{2pt}
$
  \begin{array}{l}
\\[0.05in]
[ 
\FlaTwoByTwoSingleLine 
{ A_{11} }{ A_{12} } 
{ A_{21} }{ A_{22} }, 
T_1, s_1, 
\FlaTwoByOneSingleLine{ v_1 }{ v_2}
,
\FlaOneByTwoSingleLine{W_1}{W_2} ] 
\\
~~~
  := \mbox{\sc   HQRP\_unb\_var3 }\\
~~~~~~( 
\FlaTwoByTwoSingleLine 
{ A_{11} }{ A_{12} } 
{ A_{21} }{ A_{22} }, 
T_1, s_1, 
\FlaTwoByOneSingleLine{ v_1 }{ v_2}
,
\FlaOneByTwoSingleLine{W_1}{W_2}, b ) 
\\[0.4in]
A_{22} := A_{22} - A_{21} W_2
\\[0.33in]
~
  \end{array}
$
}

\renewcommand{\initialize}{$v := \mbox{\sc ComputeWeights}( A )$}

%%% \begin{figure}[tbp]
%%% %%% \begin{algorithm}
%%% %%% $\left[ A, t \right] \becomes \URt( A )$
%%% %%% \end{algorithm}
%%% \footnotesize
%%% \FlaAlgorithmWithInit
%%% \caption{Blocked HQRP algorithm.  Note: $ W $ starts as a $ b \times
%%%   n( A ) 
%%% $ matrix.  If this is not a uniform block size used during
%%% computation, resizing may be necessary.}
%%% \label{fig:RRQR_blk}
%%% \end{figure}

%% file: 03HQRRP.tex
\input QRP_blk

\NoShow{
The fundamental problem with the blocked {\sc HQRP} algorithm is that
to cast most computation in terms of matrix-matrix multiplication, one
has to know how to swap the columns that become the column panel $
\FlaTwoByOneSingleLine{ A_{11} }{ A_{21} } $.
}
This section describes a computationally efficient technique for picking
a selection of $b$ columns from a given $n\times n$ matrix $A$ that form
good choices for the first $b$ pivots in a blocked HQRP algorithm.
Observe that this task is closely related to the problem of finding an index set $s$ of
length $b$ such that the columns in $A (:,s)$ form a good approximate basis
for the column space of $A$. Another way of expressing this problem is that we
are looking for a collection of $b$ columns whose spanning volume in $\mathbb{C}^{n}$
is close to maximal. To find the absolutely optimal choice here is a hard problem
\cite{gu1996}, but luckily, for pivoting purposes it is sufficient to find a choice
that is ``good enough.''

\subsection{Randomized pivot selection}
\label{sec:basicidea}
The strategy that we propose is simple. The idea is to perform classical QR factorization with column pivoting (QRP)
on a matrix $Y$ that is much smaller than $A$, so that performing QRP with that matrix constitutes a lower order
cost. As a bonus, it may fit in fast cache memory. This smaller matrix can be constructed by forming
random linear combinations of the rows of $ A$ as follows:
\begin{enumerate}
\item Fix an over-sampling parameter $p$. Setting $p=5$ or $p=10$ are good choices.
\item Form a random matrix $\mtx{G}$ of size $(b+p)\times n$ whose entries are drawn
independently from a normalized Gaussian distribution.
\item Form the $(b+p)\times n$ sampling matrix $\mtx{Y} = \mtx{G}\mtx{A}$.
\end{enumerate}
The sampling matrix $\mtx{Y}$ has as many columns as $\mtx{A}$, but many fewer rows.
Now execute $b$ steps of a column pivoted QR factorization
to determine an integer vector with $ b $ elements that capture how columns need to be pivoted:
\[
s = \mbox{\sc DeterminePivots}(\mtx{Y},b).
\]
In other words, the columns $Y(:,s)$ are good pivot columns for $Y$.
Our claim is that due to the way $Y$ is constructed, the columns
$A(:,s)$ are then also good choices for pivot columns of $\mtx{A}$.
This claim is supported by extensive numerical experiments,
some of which are given in Section \ref{sec:pivotquality}. There is theory supporting the claim
that these $b$ columns form a good approximate basis for the column space of $\mtx{A}$,
see, e.g.~\cite[Sec.~5.2]{2011_martinsson_randomsurvey} and \cite{2007_martinsson_PNAS,2014_martinsson_CUR},
but it has not been directly proven that they form good choices as pivots in a QR
factorization. This should not be surprising given that it is \textit{known} that even classical
column pivoting can result in poor choices \cite{1966_kahan_NLA}. Known algorithms that are provably good are
all far more complex to implement \cite{gu1996}.

Notice that there are many choices of algorithms that can be employed to determine the pivots.
For example, since high numerical accuracy is not necessary, the classical Modified Gram-Schmidt
(MGS) algorithm with column pivoting is a simple yet effective choice.

The randomized strategy described here for determining a block of pivot vectors is inspired by a
technique published in \cite[Sec.~4.1]{2006_martinsson_random1_orig} for computing a low-rank approximation
to a matrix, and later elaborated in \cite{2007_martinsson_PNAS,2011_martinsson_randomsurvey,2011_martinsson_random1}.

\begin{remark}[Choice of over-sampling parameter $p$]
The reliability of the procedures described in this section depends
on the choice of the over-sampling parameter $p$. It is well understood
how large $p$ needs to be in order to determine a high-quality
approximate basis for the column space of $\mtx{A}$ with extremely high reliability:
the choice $p=5$ is very good, $p=10$ is excellent, and $p=b$ is almost
always over-kill \cite{2011_martinsson_randomsurvey}. The pivot selection
problem is less well studied, but is more forgiving. (The choice of pivots
does not necessarily have to be     particularly optimal.) Numerical experiments indicate
that even setting $p=0$ typically
results in good choices. However, the choices $p=5$ or $p=10$ appear to be good
generic values that have resulted in excellent choices in every experiment we have run.
\end{remark}

\begin{remark}[Intuition of random projections]
To understand why the pivot columns selected by processing the small matrix $\mtx{Y}$
also form good choices for the original matrix $\mtx{A}$, it might be helpful to observe that
for a Gaussian random matrix $\mtx{G}$ of size $\ell\times n$, it is the case that for any
$x \in \mathbb{R}^{n}$, we have $\mathbb{E}\bigl[\|Gx\|^{2}\bigr] = \|x\|^{2}$, where $\mathbb{E}$
denotes expectation. Moreover, as the number of rows $\ell$ grows, the probability
distribution of $\|Gx\|$ concentrates tightly around its expected value, see, e.g.,
\cite[Sec.~2.4]{woodruff2014sketching} and the references therin. This means that for
any pair of indices $i,j \in \{1,2,\dots,n\}$ we have
$\mathbb{E}\bigl[\|Y(:,i) - Y(:,j)\|^{2}\bigr] = \|A(:,i) - A(:,j)\|^{2}$. This simple observation
does not in any way provide a proof that the randomized strategy we propose works, but might
help understand the underlying intuition.
\end{remark}

\subsection{Efficient downdating of the sampling matrix $ Y $}
\label{sec:downdateinformal}
For the QRP factorization algorithm, it is well known
that one does not need to recompute the column norms of the remainder matrix
after each step. Instead, these can be cheaply downdated, as described, e.g.,
in \cite[Ch.5, Sec.~2.1]{1998_stewart_volume1}.  In terms of asymptotic flop
counts, this observation makes the cost of pivoting become a lower order term,
and consequently both unpivoted and pivoted Householder QR algorithms have the
same leading order term $(4/3)n^3$ in their asymptotic flop\footnote{We use the
standard convention of counting one multiply and one add as one flop, regardless
of whether a complex or real operation is performed.} counts for $n\times n$ matrices.
In this section,
we describe an analogous technique for the randomized sampling strategy described
in Section \ref{sec:basicidea}. This downdating strategy was discovered by one of
the authors; a closely related technique was discovered independently and
published to arXiv by Duersch and Gu \cite{2015_blockQR_ming} in September 2015.

First observe that if the randomized sampling technique described in
Section \ref{sec:basicidea} is used in the obvious fashion, then each step
of the iteration requires the generation of a Gaussian random matrix $G$ and
a matrix-matrix multiply involving the remaining portion of $A$ in the
lower right corner to form the sampling matrix $Y$. The number of flops
required by the matrix-matrix multiplications add up to an $O(n^3)$ term
for $n\times n$ matrices.
However, it turns out to be possible to avoid computing a sampling matrix $Y$ from
scratch at every step. The idea is that if we select the randomizing matrix
$G$ in a particular way in every step beyond the first, then the corresponding
sampling matrix $Y$ can inexpensively be computed by downdating the sampling
matrix from the previous step.

To illustrate, suppose that we start with an $n\times n$ original matrix
$\mtx{A} = \mtx{A}^{(0)}$. In the first blocked step, we draw a $(b+p)\times n$
randomizing matrix $\mtx{G}^{(1)}$ and form the $(b+p)\times n$ sampling matrix
\begin{equation}
\label{eq:idea1}
\mtx{Y}^{(1)} = \mtx{G}^{(1)}\mtx{A}^{(0)}.
\end{equation}
Using the information in $\mtx{Y}^{(1)}$, we identify the $b$ pivot vectors and form
the corresponding permutation matrix $\mtx{P}^{(1)}$. Then the matrix $\mtx{Q}^{(1)}$
representing the $b$ Householder reflectors dictated by the $b$ pivot columns is formed.
Applying these transforms to the right and the left of $\mtx{A}^{(0)}$, we obtain the matrix
\begin{equation}
\label{eq:idea2}
\mtx{A}^{(1)} = \bigl(\mtx{Q}^{(1)}\bigr)^{*}\mtx{A}^{(0)}\mtx{P}^{(1)}.
\end{equation}
To select the pivots in the next step, we need to form a randomizing matrix
$\mtx{G}^{(2)}$ and a sampling matrix $\mtx{Y}^{(2)}$ that are related through
\begin{equation}
\label{eq:idea3}
\mtx{Y}^{(2)} = \mtx{G}^{(2)}\bigl(\mtx{A}^{(1)} - \mtx{R}^{(1)}\bigr),
\end{equation}
where $\mtx{R}^{(1)}$ holds the top $b$ rows of $\mtx{A}^{(1)}$ so that
$$
\mtx{A}^{(1)} - \mtx{R}^{(1)}
=
\left(\begin{array}{cc}
\mtx{A}_{11}^{(1)} & \mtx{A}_{12}^{(1)} \\
\mtx{0} & \mtx{A}_{22}^{(1)}
\end{array}\right)
-
\left(\begin{array}{cc}
\mtx{A}_{11}^{(1)} & \mtx{A}_{12}^{(1)} \\
\mtx{0} & \mtx{0}
\end{array}\right)
=
\left(\begin{array}{cc}
\mtx{0} & \mtx{0} \\
\mtx{0} & \mtx{A}_{22}^{(1)}
\end{array}\right).
$$
The key idea is now to \textit{choose} the randomizing matrix $\mtx{G}^{(2)}$
according to the formula
\begin{equation}
\label{eq:idea4}
\mtx{G}^{(2)} = \mtx{G}^{(1)}\mtx{Q}^{(1)}.
\end{equation}
Inserting (\ref{eq:idea4}) into (\ref{eq:idea3}), we now find that the sampling matrix is
\begin{multline}
\label{eq:idea5}
\mtx{Y}^{(2)} =
\mtx{G}^{(1)}\mtx{Q}^{(1)}\bigl(\mtx{A}^{(1)} - \mtx{R}^{(1)}\bigr) = \{\mbox{Use (\ref{eq:idea2})}\} = \\
\mtx{G}^{(1)}\mtx{A}^{(0)}\mtx{P}^{(1)} - \mtx{G}^{(1)}\mtx{Q}^{(1)}\mtx{R}^{(1)} = \{\mbox{Use (\ref{eq:idea1})}\} =
\mtx{Y}^{(1)}\mtx{P}^{(1)} - \mtx{G}^{(1)}\mtx{Q}^{(1)}\mtx{R}^{(1)}.
\end{multline}
Evaluating formula (\ref{eq:idea5}) is inexpensive since the first term is
a permutation of the columns of the sampling matrix $\mtx{Y}^{(1)}$ and the
second term is a product of thin matrices (recall that $\mtx{Q}^{(1)}$
is a product of $b$ Householder reflectors).

\begin{remark}
Since the probability distribution for Gaussian random matrices is invariant under
unitary maps, the formula (\ref{eq:idea4}) appears quite safe. After all,
$\mtx{G}^{(1)}$ is Gaussian, and $\mtx{Q}^{(1)}$
is just a sequence of reflections, so it might be tempting to conclude that the new
randomizing matrix must be Gaussian too. However, the matrix $\mtx{Q}^{(1)}$
unfortunately depends on the draw of $\mtx{G}^{(1)}$, so this argument does
not work. Nevertheless, the dependence of $\mtx{Q}^{(1)}$ on $\mtx{G}^{(1)}$
is very subtle since this $\mtx{Q}^{(1)}$ is dictated primarily by the directions of the good pivot columns.
Extensive practical experiments (see, e.g., Section \ref{sec:pivotquality}) indicate that the
pivoting strategy described in this section based on downdating is just as good as the one that
uses ``pure'' Gaussian matrices that was described in Section \ref{sec:basicidea}.
\end{remark}

\subsection{Detailed description of the downdating procedure}
\label{sec:downdating}

Having described the downdating procedure informally in Section \ref{sec:downdateinformal},
we in this section provide a detailed description using the notation for \HQR~that we
used in Section \ref{sec:HPQRRP}. First, let us assume that one iteration
of the blocked algorithm has completed, so that, at the bottom of the loop body,
the matrix $ A $ contains
\[
\FlaTwoByTwoSingleLine{ A_{11} }{ A_{12} }{ A_{21} }{ A_{22} } =
\FlaTwoByTwoSingleLine{ { U \backslash R }_{11} }{ R_{12} }{ U_{21} }{ \widehat A_{22} - U_{21} W_{12} } .
\]
Here $ \widehat{ A} $ denotes the original contents of $ A P_1 $, where $ P_1 $ captures how columns have been swapped so far.
Hence, cf.~(\ref{eq:idea2}),
\[
\begin{array}[t]{c}
\underbrace{
\left(  I - \FlaTwoByOneSingleLine{U_{11}}{U_{21}}
T_{11}^{-H} \FlaTwoByOneSingleLine{U_{11}}{U_{21}}^H \right)
} \\
\RvdG{\left( \begin{array}{c | c}
Q_1 & Q_2
\end{array}
\right)}
\end{array}
\FlaTwoByTwoSingleLine{ \widehat A_{11} }{ \widehat A_{12} }{ \widehat A_{21} }{ \widehat A_{22} } P_1 =
\FlaTwoByTwoSingleLine{ R_{11} }{ R_{12} }{ 0 }{ A_{22} }
\]

\RvdG{
Now, let $ \widetilde G_2 $ be the next sampling matrix and
$ \widetilde Y_2 = \widetilde G_2  A_{22} $.
In order to show how this new sampling matrix can be computed by
downdating the last sampling matrix, consider that
\[
\FlaOneByTwoSingleLine{\widetilde Y_1}{\widetilde Y_2} =
\FlaOneByTwoSingleLine{\widetilde G_1}{\widetilde G_2}
\FlaTwoByTwoSingleLine{ 0 }{ 0 }{ 0 }{ A_{22} }
\]
for some matrix $ \widetilde G_1 $ \RvdG{and that}
}%
\NoShow{Letting
$$
\FlaOneByTwoSingleLine{ Q_1}{Q_2} = \left(  I - \FlaTwoByOneSingleLine{U_{11}}{U_{21}}
T_{11}^{-1} \FlaTwoByOneSingleLine{U_{11}}{U_{21}}^H \right)
$$
we find that
}
\begin{multline}
\lefteqn{
\FlaOneByTwoSingleLine{\widetilde G_1}{\widetilde G_2}\FlaTwoByTwoSingleLine{ 0 }{ 0 }{ 0 }{ A_{22} }
=
\FlaOneByTwoSingleLine{\widetilde G_1}{\widetilde G_2}
\left(
\FlaTwoByTwoSingleLine{ R_{11}  }{ R_{12} }{ 0 }{ A_{22} }
-
\FlaTwoByTwoSingleLine{ R_{11}  }{ R_{12} }{ 0 }{ 0 }
\right)}
\\
\NoShow{
=
\FlaOneByTwoSingleLine{\widetilde G_1}{\widetilde G_2}
\begin{array}[t]{c}
\underbrace{
\FlaOneByTwoSingleLine{ Q_1}{Q_2}^H \FlaOneByTwoSingleLine{ Q_1}{Q_2}}
\\
I
\end{array}
\left(
\FlaTwoByTwoSingleLine{ R_{11}  }{ R_{12} }{ 0 }{ A_{22} }
-
\FlaTwoByTwoSingleLine{ {R}_{11}  }{ R_{12} }{ 0 }{ 0 }
\right)
\\
}
=
\FlaOneByTwoSingleLine{\widetilde G_1}{\widetilde G_2}
\FlaTwoByTwoSingleLine{ R_{11}  }{ R_{12} }{ 0 }{ A_{22} }
-
\FlaOneByTwoSingleLine{\widetilde G_1}{\widetilde G_2}
\FlaTwoByTwoSingleLine{ {R}_{11}  }{ R_{12} }{ 0 }{ 0 }
\\
=
\FlaOneByTwoSingleLine{\widetilde G_1}{\widetilde G_2}
\FlaOneByTwoSingleLine{ Q_1}{Q_2} \widehat A P_1 -
\FlaOneByTwoSingleLine{\widetilde G_1}{\widetilde G_2}
\FlaOneByTwoSingleLine{ Q_1}{Q_2}
\FlaOneByTwoSingleLine{ Q_1}{Q_2}^H
\FlaTwoByTwoSingleLine{ {R}_{11}  }{ R_{12} }{ 0 }{ 0 }
\label{eq:swiss}
\end{multline}
The choice of randomizing matrix analogous to (\ref{eq:idea4}) is now
\begin{equation}
\label{eq:swisschoice}
\FlaOneByTwoSingleLine{\widetilde G_1}{\widetilde G_2} =
\FlaOneByTwoSingleLine{G_1}{G_2}
\FlaOneByTwoSingleLine{Q_1}{Q_2}.
\end{equation}
Inserting the choice (\ref{eq:swisschoice}) into (\ref{eq:swiss}), we obtain
\begin{eqnarray*}
\FlaOneByTwoSingleLine{\widetilde Y_1}{\widetilde Y_2}
=
\begin{array}[t]{c}
\underbrace{
\FlaOneByTwoSingleLine{G_1}{G_2}
\widehat A
}
\\
\FlaOneByTwoSingleLine{Y_1}{Y_2}
\end{array}
 P_1 -
\FlaOneByTwoSingleLine{G_1}{G_2}
\left(  I - \FlaTwoByOneSingleLine{U_{11}}{U_{21}}
T_{11}^{-H} \FlaTwoByOneSingleLine{U_{11}}{U_{21}}^H \right)^H
\FlaTwoByTwoSingleLine{ {R}_{11}  }{ R_{12} }{ 0 }{ 0 } .
\end{eqnarray*}
Letting
$ \FlaOneByTwoSingleLine{\overline Y_1}{\overline Y_2} =
\FlaOneByTwoSingleLine{Y_1}{Y_2} P_1 $
we conclude that
\begin{eqnarray*}
\widetilde Y_2 &=&
\overline Y_2 -
\FlaOneByTwoSingleLine{G_1}{G_2}
\left(  I - \FlaTwoByOneSingleLine{U_{11}}{U_{21}}
T_{11}^{-1} \FlaTwoByOneSingleLine{U_{11}}{U_{21}}^H \right)
\FlaTwoByOneSingleLine{ R_{12} }{ 0 }\\
&=&
\overline Y_2 -
\FlaOneByTwoSingleLine{G_1}{G_2}
\left(  \FlaTwoByOneSingleLine{ R_{12} }{ 0 }
 - \FlaTwoByOneSingleLine{U_{11}}{U_{21}}
T_{11}^{-1} {U_{11}}^H R_{12} \right) \\
&=&
\overline Y_2 -
\left(  G_1
 - ( G_1 U_{11} + G_2 {U_{21}} )
T_{11}^{-1} {U_{11}}^H  \right) { R_{12} }.
\end{eqnarray*}
which can then be used to downdate the sampling
matrix $ Y $.

\subsection{The blocked algorithm}

In Fig.~\ref{fig:HQRRPblk}, we give the blocked algorithm that results when
the randomized pivot selection strategy described in Section \ref{sec:basicidea}
is combined with the downdating techniques described in Section \ref{sec:downdating}.
In that figure, there is a call to a function ``HRQP\_{\small{UNB}}'' which is
an unblocked \HQR\ algorithm with column pivoting. The purpose of this call is 
to factor the current column panel so that the diagonal elements within blocks 
on the diagonal of $R$ are ordered from largest to smallest in magnitude.
Moreover, the call to ``UpdatePivotInfo'' takes the pivoting that occurred within the current panel (to ensure strictly decreasing diagonal elements within the current diagonal block of $R_{11}$) and merges this with the pivot information that occurred when determining the columns to be moved into that current panel.

\subsection{Asymptotic cost analysis}

\input 031cost_new

%% file: QRP_blk.tex
\setlength{\arraycolsep}{2pt}
\resetsteps      % Reset all the commands to create a blank worksheet  

% Define the operation to be computed

\renewcommand{\routinename}{ \left[ A, T, s \right] := \mbox{\sc HQRP\_randomized\_blk}( A, T, s, b, p ) }

% Step 3: Loop-guard 

\renewcommand{\initialize}{
$
\begin{array}{@{}l}
 G := \mbox{\sc rand\_iid}( b+p, n(A) )  \\
 Y := G A \\[-0.05in]
\end{array}
$
}

\renewcommand{\guard}{
  m( A_{TL} ) < m( A )
}

% Step 4: Define Initialize 

\renewcommand{\partitionings}{
  $
  A \rightarrow
  \FlaTwoByTwo{A_{TL}}{A_{TR}}
              {A_{BL}}{A_{BR}}
  $
,
  $
  T \rightarrow
  \FlaTwoByOne{T_{T}}
              {T_{B}}
  $
,
  $
  s \rightarrow
  \FlaTwoByOne{s_{T}}
              {s_{B}}
  $
,
  $
  Y \rightarrow
  \FlaOneByTwo{Y_L}{Y_R}
  $
}

\renewcommand{\partitionsizes}{
$ A_{TL} $ is $ 0 \times 0 $,
$ T_{T} $ has $ 0 $ rows,
$ s_{T} $ has $ 0 $ rows,
$ Y_L $ has $ 0 $ columns
}

% Step 5a: Repartition the operands 

\renewcommand{\blocksize}{b \rightarrow \min(b,n(A_{BR}))}

\renewcommand{\repartitionings}{
$  \FlaTwoByTwo{A_{TL}}{A_{TR}}
              {A_{BL}}{A_{BR}}
  \rightarrow
  \FlaThreeByThreeBR{A_{00}}{A_{01}}{A_{02}}
                    {A_{10}}{A_{11}}{A_{12}}
                    {A_{20}}{A_{21}}{A_{22}}
$,
$  \FlaTwoByOne{ T_T }
              { T_B }
\rightarrow
  \FlaThreeByOneB{T_0}
                 {T_1}
                 {T_2}
$
,\\
$  \FlaTwoByOne{ s_T }
              { s_B }
\rightarrow
  \FlaThreeByOneB{s_0}
                 {s_1}
                 {s_2}
$
,
$  \FlaOneByTwo{Y_L}{Y_R}
\rightarrow  \FlaOneByThreeR{Y_0}{Y_1}{Y_2}
$
}

\renewcommand{\repartitionsizes}{
  $ A_{11} $ is $ b \times b $,
  $ T_1 $ has $ b $ rows,
  $ s_1 $ has $ b $ rows,
  $ Y_1 $ has $b$ columns}

% Step 5b: Move the double lines 

\renewcommand{\moveboundaries}{
$  \FlaTwoByTwo{A_{TL}}{A_{TR}}
              {A_{BL}}{A_{BR}}
  \leftarrow
  \FlaThreeByThreeTL{A_{00}}{A_{01}}{A_{02}}
                    {A_{10}}{A_{11}}{A_{12}}
                    {A_{20}}{A_{21}}{A_{22}}
$,
$  \FlaTwoByOne{ T_T }
               { T_B }
\leftarrow
  \FlaThreeByOneT{T_0}
                 {T_1}
                 {T_2}
$
, \\
$  \FlaTwoByOne{ s_T }
               { s_B }
\leftarrow
  \FlaThreeByOneT{s_0}
                 {s_1}
                 {s_2}
$
, 
$  \FlaOneByTwo{Y_L}{Y_R}
\leftarrow  \FlaOneByThreeL{Y_0}{Y_1}{Y_2}
$
}

% Step 8: Insert the updates required to change the 
%         state from that given in Step 6 to that given in Step 7
% Note: The below needs editing!!!

\renewcommand{\update}{
$
  \begin{array}{l}
      s_1  := \mbox{\sc DeterminePivots}( 
        \FlaOneByTwoSingleLine{ Y_1 }{ Y_2 }, b ) \\
\left( \begin{array}{c | c}
A_{01} & A_{02} \\ \whline
A_{11} & A_{12} \\ \hline
A_{21} & A_{22}
\end{array} \right)
:= \mbox{\sc SwapCols}( s_1, 
\left( \begin{array}{c | c}
A_{01} & A_{02} \\ \whline
A_{11} & A_{12} \\ \hline
A_{21} & A_{22}
\end{array} \right)
) \\
\rowcolor[gray]{0.88}
\left[
\FlaTwoByOneSingleLine{ A_{11} }{ A_{21} }, T_1\RvdG{\colorbox{white}{$, s_1'$}}
\right]
:= 
\mbox{\sc HQR\RvdG{\colorbox{white}{\sc P}}\_unb}( 
\FlaTwoByOneSingleLine{ A_{11} }{ A_{21} } , T_1
)
\\
\rowcolor[gray]{1}
\RvdG{
A_{01} := \mbox{\sc SwapCols}( s_1', A_{01} )
}
\\
\rowcolor[gray]{1}
\RvdG{
s_1 := \mbox{\sc UpdatePivotInfo}( s_1', s_1 )
}
\\
\rowcolor[gray]{0.88} 
W_{12} := 
T_{1}^{-H} ( U_{11}^H A_{12} + U_{21}^H A_{22} ) \\
\rowcolor[gray]{0.88}
\FlaTwoByOneSingleLine
  { A_{12} }
  { A_{22} }
:= 
\FlaTwoByOneSingleLine
  { A_{12} - U_{11} W_{12} }
  { A_{22} - U_{21} W_{12} }
%%% \FlaTwoByOneSingleLine{ A_{12} }{ A_{22} } 
%%% :=
%%% \mbox{HQR\_Apply}( T_{1, 
%%% \FlaTwoByOneSingleLine{ A_{11} }{ A_{21} },
%%% \FlaTwoByOneSingleLine{ A_{12} }{ A_{22} } )
\\
\FlaOneByTwoSingleLine{ Y_1 }{ Y_2 }
:=
\mbox{\sc SwapCols}( s_1, \FlaOneByTwoSingleLine{ Y_1 }{ Y_2 } )\\
Y_2
:=
Y_2 - 
\left(  G_1 
 - ( G_1 U_{11} + G_2 {U_{21}} )
T_{11}^{-1} {U_{11}}^H  \right) { R_{12} }.
  \end{array}
$
}

\begin{figure}[tbp]
\FlaAlgorithmWithInit
\caption{Blocked Householder transformation based QR factorization
  with column pivoting based on randomization.  In this algorithm, $
  U_{11} $ equals the unit lower triangiular matrix stored below the
  diagonal of $ A_{11} $,  $ U_{21} = A_{21} $, and $ R_{12} = A_{12}
  $.  The steps highlighted in gray constitute the blocked QR
  factorization without column pivoting from Fig.~\ref{fig:HQR_unb_blk}.
}
\label{fig:HQRRPblk}
\end{figure}

%% file: 031cost_new.tex
In analyzing the asymptotic complexity of the method, we consider a
matrix of size $m\times n$,  with $m \geq n$. We assume that the block size
$b$ and the over-sampling parameter $p$ are kept fixed as $m$ and $n$ grow.
We first note that all steps in Fig.~\ref{fig:HQRRPblk}
highlighted in grey are part of the blocked HQR
algorithm, which is known to have an asymptotic cost of $2mn^{2} - 2/3 n^3 $ flops.
This leaves us to discuss the overhead related to the other
operations.
\begin{itemize}
\item
{$ G := \mbox{\sc rand\_iid}( b + p, m) $:}
Cost: ignored.
\item
{$ Y := G A $:}
Cost: $O( ( b + p ) mn) $ flops.
\item
{$ s_1  := \mbox{\sc DeterminePivots}(
        \FlaOneByTwoSingleLine{ Y_1 }{ Y_2 }, b )  $:}
Cost: $ O( b (b + p ) ( n - kb ) )$ flops
\RvdG{during the $ k $th iteration of the blocked algorithm, for a total of
$ O( (b+p)n^2)$ flops}.  (Recall that the factorization of this matrix
can stop after the first $ b $ columns have been identified.)
%% \item
%% $ \left( \begin{array}{c | c}
%% A_{01} & A_{02} \\ \whline
%% A_{11} & A_{12} \\ \hline
%% A_{21} & A_{22}
%% \end{array} \right)
%% := \mbox{\sc SwapCols}( s_1,
%% \left( \begin{array}{c | c}
%% A_{01} & A_{02} \\ \whline
%% A_{11} & A_{12} \\ \hline
%% A_{21} & A_{22}
%% \end{array} \right)
%% ) $\\
%% Cost: ignored.
\item
{$ \cdots
:=
\mbox{\sc SwapCols}( s_1, \cdots  ) $:}
Cost: ignored.
\item
{$ Y_2
:=
Y_2 -
\left(  G_1
 - ( G_1 U_{11} + G_2 {U_{21}} )
T_{11}^{-1} {U_{11}}^H  \right) { R_{12} } $:}
Aggregate cost: $ O( (b+p) n^2 ) $.
\end{itemize}
Thus,
the overhead is
$ O( (b+p)  (n^2 + mn))$ flops and the total cost is
\[
2mn^{2} - 2/3 n^3
+
O( (b+p)  (n^2 + mn))
\mbox{ flops},
\]
which, asymptotically, approaches the same $2mn^2 - 2/3 n^3 $ cost as unpivoted \HQR.

%% file: 04experiments.tex
This section describes the results from two sets of experiments.
Section \ref{sec:performance} compares the computational speed of
the proposed scheme to existing state-of-the-art methods for computing
column pivoted QR factorizations.
Section \ref{sec:pivotquality} investigates how well the proposed
randomized technique works at selecting pivot columns.
Specifically,
we investigate how well the rank-$k$ truncated QR factorization
approximates the original matrix and compare the results to those
obtained by classical column pivoting.

\subsection{Performance experiments}
\label{sec:performance}

We have implemented the proposed HQRRP algorithm using the {\tt libflame}~\cite{CiSE09,libflame_ref}
library that allows our implementations to closely resemble the algorithms as presented in the paper.

\vspace{0.05in}
\noindent
{\bf Platform details.}
All experiments reported in this article were performed on an
Intel Xeon E5-2695 v3 (Haswell) processor (2.3 GHz), with 14 cores.
In order to be able to show scalability results, the clock speed was
throttled at 2.3 GHz, turning off so-called turbo boost.
Each core can execute 16  double precision  floating point operations
per cycle, meaning that the peak performance of one core is
36.8 GFLOPS (billions of floating point operations per second).
For comparison, on a single core, {\tt dgemm} achieves around 33.6 GFLOPS.
Other details of interest
include
that the OS used was Linux (Version 2.6.32-504.el6.x86\_64),
the code was compiled with gcc (Version 4.4.7),
{\tt dgeqrf} and {\tt dgeqp3} were taken from LAPACK (Release 3.4.0),
and the implementations were linked to BLAS from Intel's MKL library
(Version 11.2.3).
Our implementations were coded with {\tt libflame} (Release 11104).

\vspace{0.05in}
\noindent
{\bf Implementations.}
We report performance for four implementations:
\begin{description}
\item
{\bf {\tt dgeqrf}.}  The implementation of blocked \HQR\ that is part
of the {\tt netlib} implemenation of LAPACK, modified so that the
block size can be controlled.
\item
{\bf {\tt dgeqp3}.}  The implementation of blocked \HQRP\ that is part of the {\tt netlib} implemenation of LAPACK, based on~\cite{QRP:SIAM}, modified so that the
block size can be controlled.
\item
{\bf {\tt HQRRPbasic}.}  Our implementation of \HQRRP\ that computes new
matrices $ G $ and $ Y $ in every iteration.   This implementation
deviates from the algorithm in Fig.~\ref{fig:HQRRPblk} in that it
also incorporates additional column pivoting within the call to
{\sc HQRRP\_unb}.
\item
{\bf {\tt HQRRP}.}  The implementation of \HQRRP\ that
downdates $ Y $. (It also includes pivoting within {\sc HQRRP\_unb}).
\item
{\bf {\tt dgeqpx}.}  An implementation of \HQRP\ with window pivoting \cite{RRQRblk1,RRQRblk2},
briefly mentioned in the introduction.
This algorithm consists of two stages:
The first stage is a QR with window pivoting. An incremental condition
estimator is employed to accept/reject the columns within a window.
The size of the window is about about twice the block size used by the algorithm to
maintain locality.
If all the columns in the window are unacceptable, all of them are rejected and
fresh ones are brought into the window.
The second stage is a postprocessing stage to validate the rank,
that is, to check that all good columns are in the front, and all
the bad columns are in the rear. This step is required because the
window pivoting could fail to reveal the rank due to its
shortsightedness (having only checked the window and having employed a cheap to compute condition
estimator.)  Sometimes, some columns must be moved between $ R_{11} $
that has been computed
and matrices $ R_{12} $ and $ R_{22} $,
and then retriangularization must be performed
with Givens rotations.
\end{description}
In all cases, we used algorithmic block sizes of $ b = 64$ and $ 128$.
While likely not optimal for all problem sizes, these blocks sizes yield near best
performance and, regardless, it allows us to easily compare and
contrast the performance of the different implementations.

\begin{figure}
\setlength{\unitlength}{1mm}
\begin{picture}(130,130)
\put(04,04){\includegraphics[height=115mm]{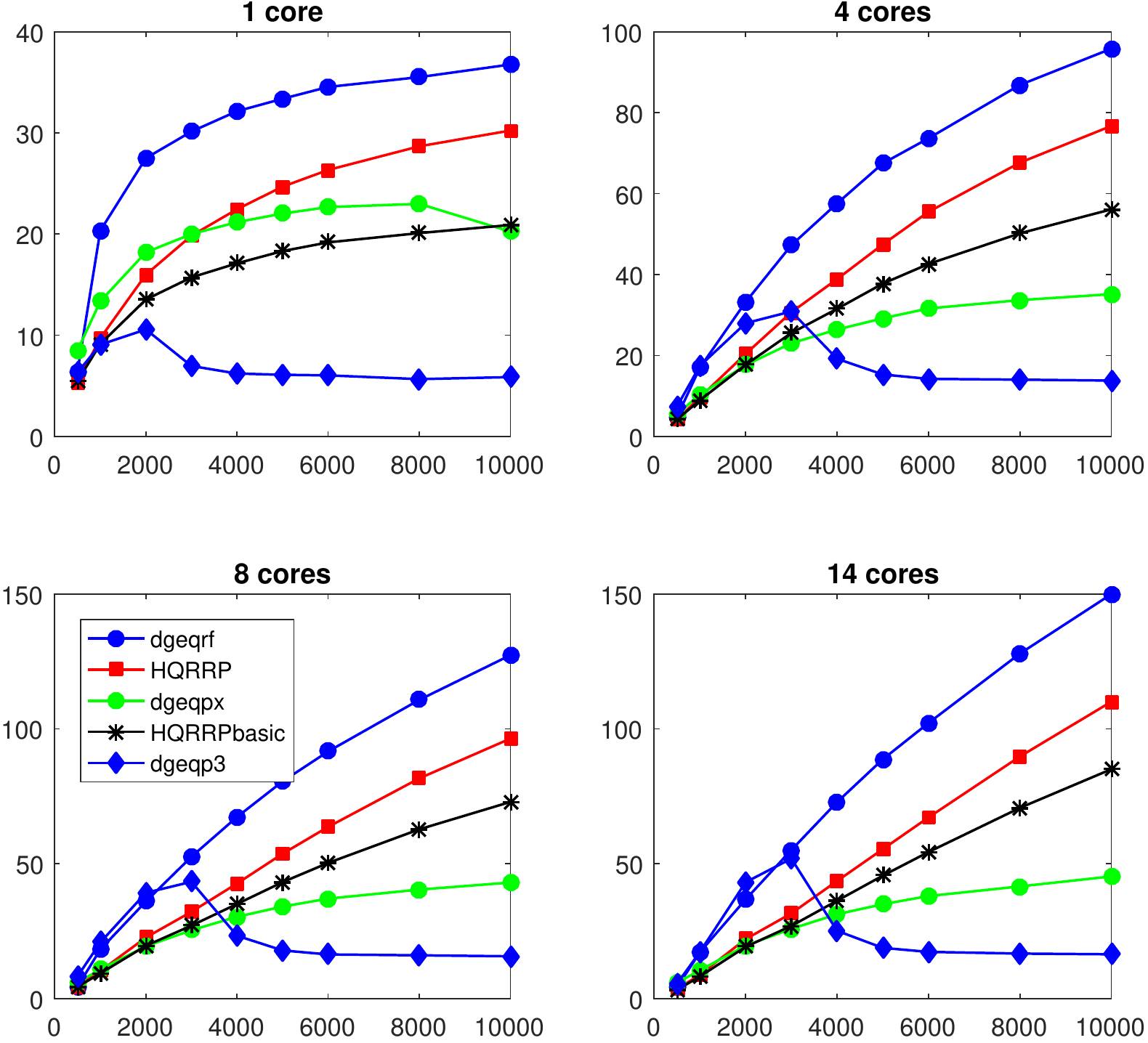}}
\put(00,27){\rotatebox{90}{GFlops}}
\put(37,00){$n$}
\put(103,00){$n$}
%\put(29,61){8 cores}
%\put(94,61){14 cores}
%\put(29,125){1 core}
%\put(94,125){4 cores}
\end{picture}
\caption{Computational speed (in standardized Gigaflops) of the proposed randomized
algorithm \HQRRP~for computing a column pivoted QR-factorization of a matrix of size
$n\times n$. The four graphs show results from test runs on 1, 4, 8, and 14 cores on
an Intel Xeon E5-2695 v3. Observe that the scales on the vertical axes are different
in the four graphs. For comparison, the graphs also show the times for competing
algorithms (\texttt{dgeqp3} and \texttt{dgeqpx}), and for unpivoted QR (\texttt{dgeqrf}),
see Section \ref{sec:performance} for details. The proposed algorithm \HQRRP~is
faster than competing algorithms, with the gap growing as more cores are added.
The block size for all algorithms was set to $b=64$. Figure \ref{fig:speedup} shows
the relative speeds in detail.}
\label{fig:performance}
\end{figure}
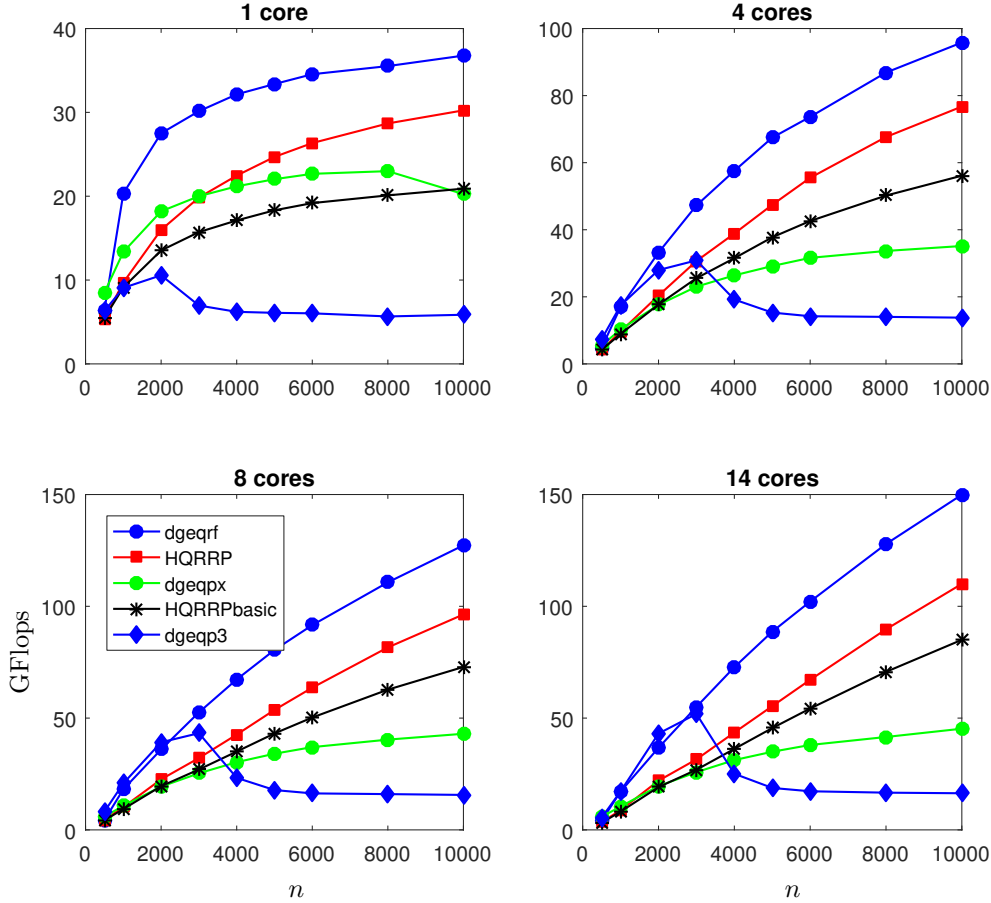

\begin{figure}
\setlength{\unitlength}{1mm}
\begin{picture}(133,66)
\put(04,04){\includegraphics[height=60mm]{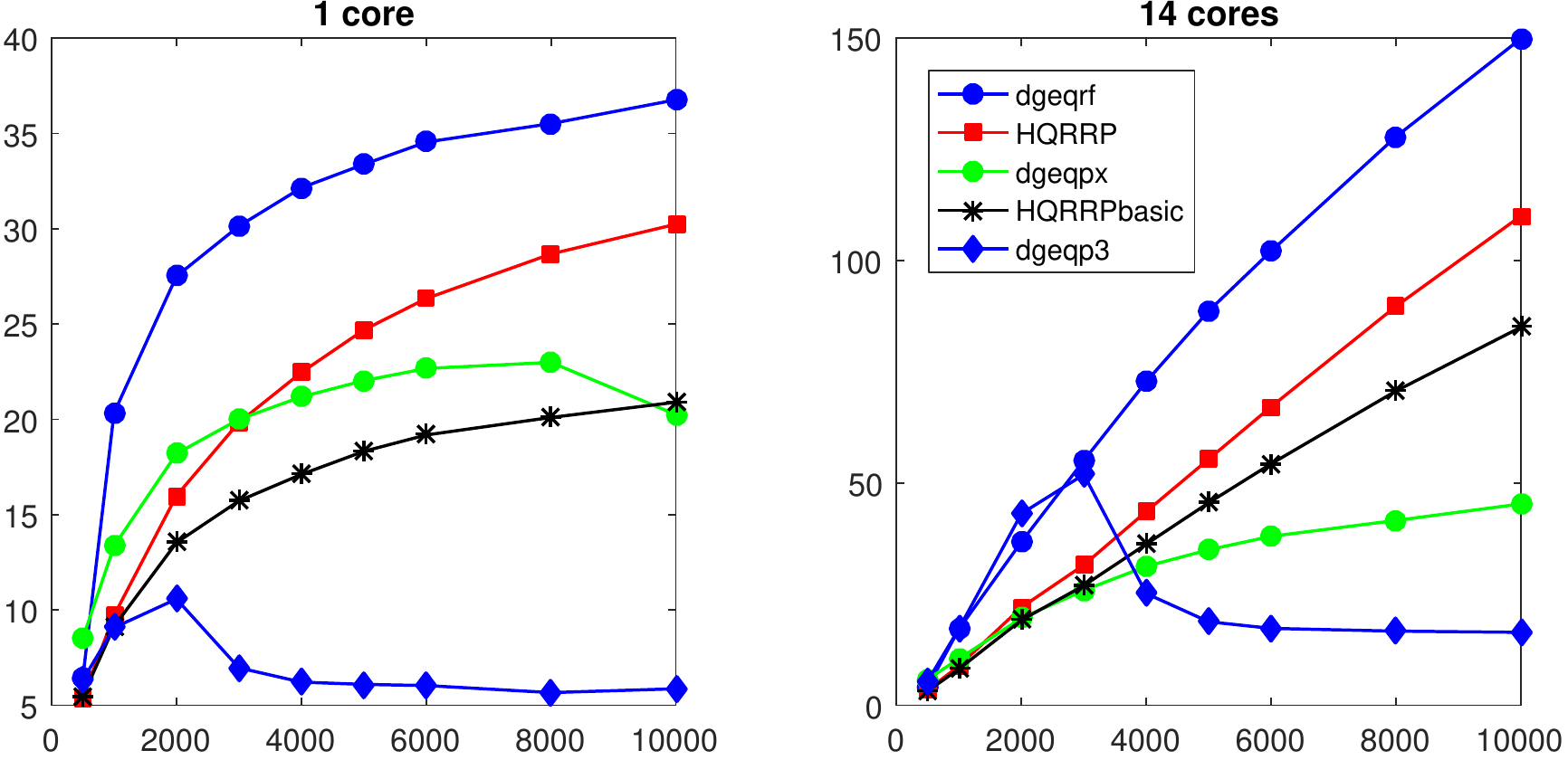}}
%\put(69,04){\includegraphics[height=55mm]{Pics/fig_speed_14c_b=128.pdf}}
\put(00,27){\rotatebox{90}{GFlops}}
\put(35,00){$n$}
\put(101,00){$n$}
%\put(29,61){1 core}
%\put(94,61){14 cores}
\end{picture}
\caption{Computational speeds (in standardized Gigaflops) for the same problem
as that shown in Figure \ref{fig:performance}. In this figure, results are shown
for a blocksize of $b=128$, in contrast with the blocksize $b=64$ that was used
in Figure \ref{fig:performance}.}
\label{fig:performance_b128}
\end{figure}

\noindent
{\bf Results.}
\RvdG{As is customary in these kinds of performance comparisons, we compute the achieved performance as%\footnote{Should it be $10^{-9}$ or $2^{-30} = 10^{-9}\times(1/1.0737\dots)$?}
\[
\frac{4/3 n^3}{\mbox{time (in sec.)}} \times 10^{-9} \mbox{ GFLOPS}.
\]
Thus, even for the implementations that perform more computations, we
only count the floating point operations performed by an unblocked
\HQR\ without pivoting.}

Figure~\ref{fig:performance} reports performance on 1, 4, 8, and 14 cores.
We see that {\tt HQRRP} handily outperforms {\tt dgeqp3} and to a lesser
degree also outperforms {\tt dgeqpx}. Moreover, the asymptotic performance
of {\tt HQRRP} appears to approach that of {\tt dgeqrf}, in particular for
the single core case. We also see that while the relative performances of
all 5 methods remain qualitatively the same across the four graphs, it is
clear that as the number of cores grows, the speed advantage of
{\tt HQRRP} over {\tt dgeqp3} becomes even further pronounced, cf.~Figure \ref{fig:speedup}.

Figure \ref{fig:performance} also shows that while the absolute speed of
all five algorithms studied improves as the number of cores grows, all
five fall substantially short of ideal scaling (in that the number of
Gigaflops \textit{per core} falls as the number of cores grows). This
observation underscores the need for further research in this area.

In order to investigate the effect of the block size on the computational
speed, we reran the experiments shown in Figure \ref{fig:performance} with
a block size of $b=128$ instead of $b=64$, with the results shown in
Figure \ref{fig:performance_b128}. We see that the choice $b=64$ tends to lead to slightly
faster execution, but the key take-away from this comparison is that the
speed is relatively insensitive to the precise choice of block size (within
reason, of course).

\subsection{Quality experiments}
\label{sec:pivotquality}
%In order to evaluate how the randomized pivot selection strategy
%proposed here fares against classical column pivoting,
In this section, we describe the results of numerical experiments that were
conducted to compare the quality of the pivot choices made by our randomized
algorithm \HQRRP~to those resulting from classical column pivoting. Specifically,
we compared how well partial factorizations reveal the numerical ranks of four
different test matrices:
%\begin{description}
\begin{itemize}
\item \textit{Matrix 1 (fast decay):} This is an $n\times n$ matrix of the form
$\mtx{A} = \mtx{U}\mtx{D}\mtx{V}^{*}$ where $\mtx{U}$ and $\mtx{V}$ are randomly
drawn matrices with orthonormal columns (obtained by performing \RvdG{QR} on a random
Gaussian matrix), and where $\mtx{D}$ is diagonal with entries
given by $d_{j} = \beta^{(j-1)/(n-1)}$ with $\beta = 10^{-5}$.
\item \textit{Matrix 2 (S shaped decay):} This matrix is built in the
same manner as ``Matrix 1'', but now the diagonal entries
of $\mtx{D}$ are chosen to first hover around 1, then decay rapidly, and then
level out at $10^{-6}$, as shown in Figure \ref{fig:errors_S} (black line).
\item \textit{Matrix 3 (Single Layer BIE):} This matrix is the result of discretizing
a Boundary Integral Equation (BIE) defined on a smooth closed curve in the plane. To be precise,
we discretized the so called ``single layer'' operator associated with the Laplace equation
using a $6^{\rm th}$ order quadrature rule designed by Alpert \cite{alpert1999hybrid}.
This is a well-known ill-conditioned operator for which column
pivoting is essential in order to stably solve the corresponding linear system.
\item \textit{Matrix 4 (Kahan):} This is a variation of the ``Kahan counter-example'' \cite{1966_kahan_NLA} which
is designed to trip up classical column pivoting.
The matrix is formed as $\mtx{A} = \mtx{S}\mtx{K}$ where:
\begin{equation*}
\mtx{S} = \left[\begin{matrix}
1 & 0 & 0 & 0 & \cdots \\
0 & \zeta & 0 & 0 & \cdots \\
0 & 0 & \zeta^2 & 0 & \cdots \\
0 & 0 & 0 & \zeta^3 & \cdots \\
\vdots & \vdots & \vdots & \vdots & \ddots \\
\end{matrix}\right] \quad \mbox{and} \quad
\mtx{K} = \begin{bmatrix}
1 & -\phi & -\phi & -\phi & \cdots \\
0 & 1 & -\phi & -\phi & \cdots \\
0 & 0 & 1 & -\phi & \cdots \\
0 & 0 & 0 & 1 & \cdots \\
\vdots & \vdots & \vdots & \vdots & \ddots \\
\end{bmatrix}
\end{equation*}
for some positive parameters $\zeta$ and $\phi$ such that $\zeta^2+\phi^2=1$.
In our experiments, we chose $\zeta = 0.99999$.
\end{itemize}

For each test matrix, we computed QR factorizations
\begin{equation}
\label{eq:AP=QRnum}
\mtx{A}\mtx{P} = \mtx{Q}\mtx{R}
\end{equation}
using three different techniques:
\begin{itemize}
\item \texttt{HQRP}: The standard QR factorization \texttt{qr} built in to Matlab R2015a.
\item \texttt{HQRRPbasic}: The randomized algorithm described in Figure \ref{fig:HQRRPblk}, but
without the updating strategy for computing the sample matrix $\mtx{Y}$.
\item \texttt{HQRRP:} The randomized algorithm described in Figure \ref{fig:HQRRPblk}.
\end{itemize}
Our implementations of both \texttt{HQRRPbasic} and \texttt{HQRRP} deviate from what is shown
in Figure \ref{fig:HQRRPblk} in that they also incorporate column pivoting within the call to
{\sc HQRRP\_unb}.
In all experiments, we used test matrices of size $4\,000\times 4\,000$, a block size of $b=100$, and
an over-sampling parameter of $p=5$.

\begin{figure}
\includegraphics[width=\textwidth]{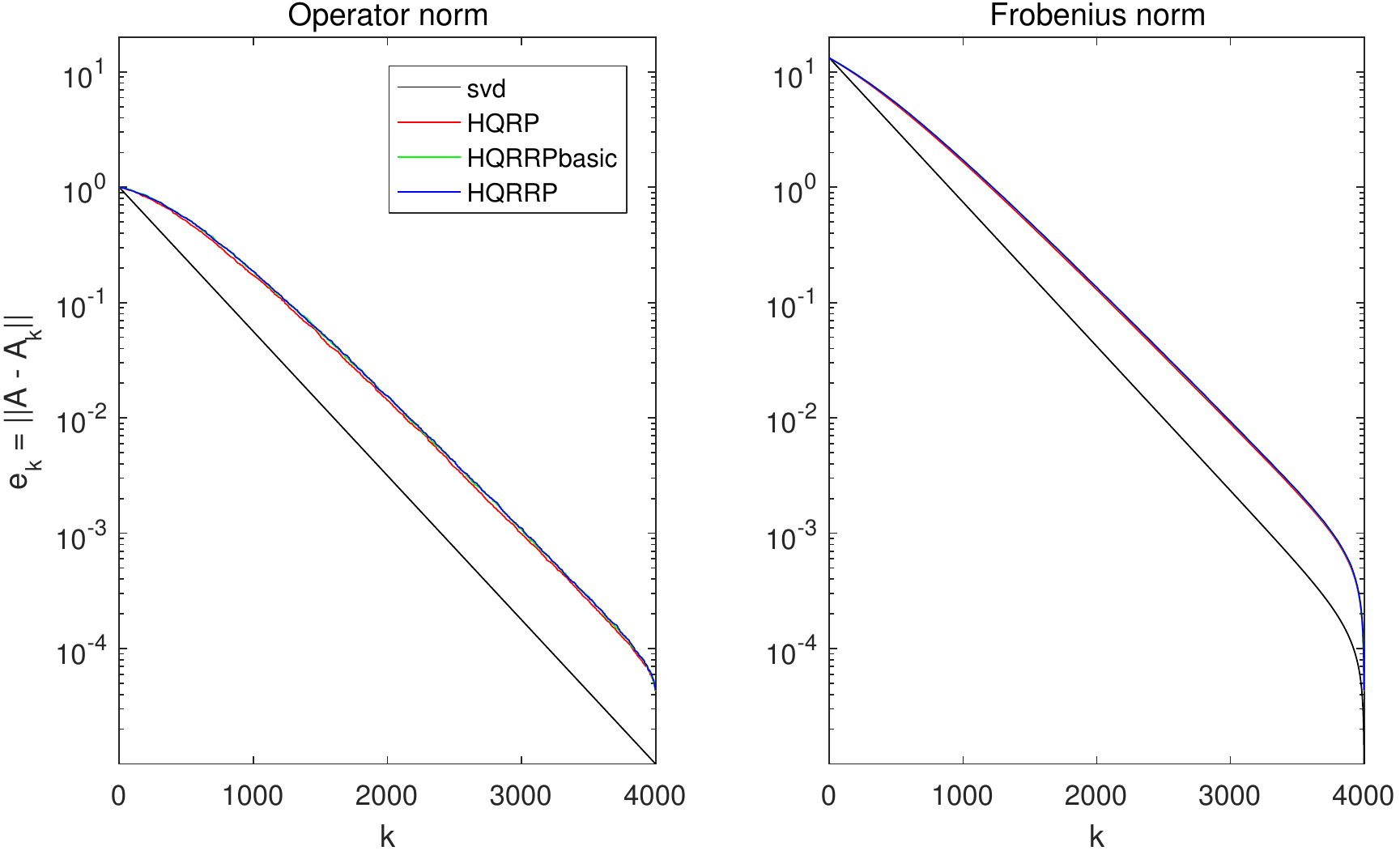}
\caption{Residual errors $e_{k}$ for ``Matrix 1'' as a function of the truncation rank $k$,
cf.~(\ref{eq:def_ek}). The red line shows the results from traditional
column pivoting, while the green and blue lines refer to the randomized methods
we propose. The black line indicates the theoretically minimal errors resulting
from a rank-$k$ partial singular value decomposition.}
\label{fig:errors_fast}
\end{figure}

\begin{figure}
\includegraphics[width=\textwidth]{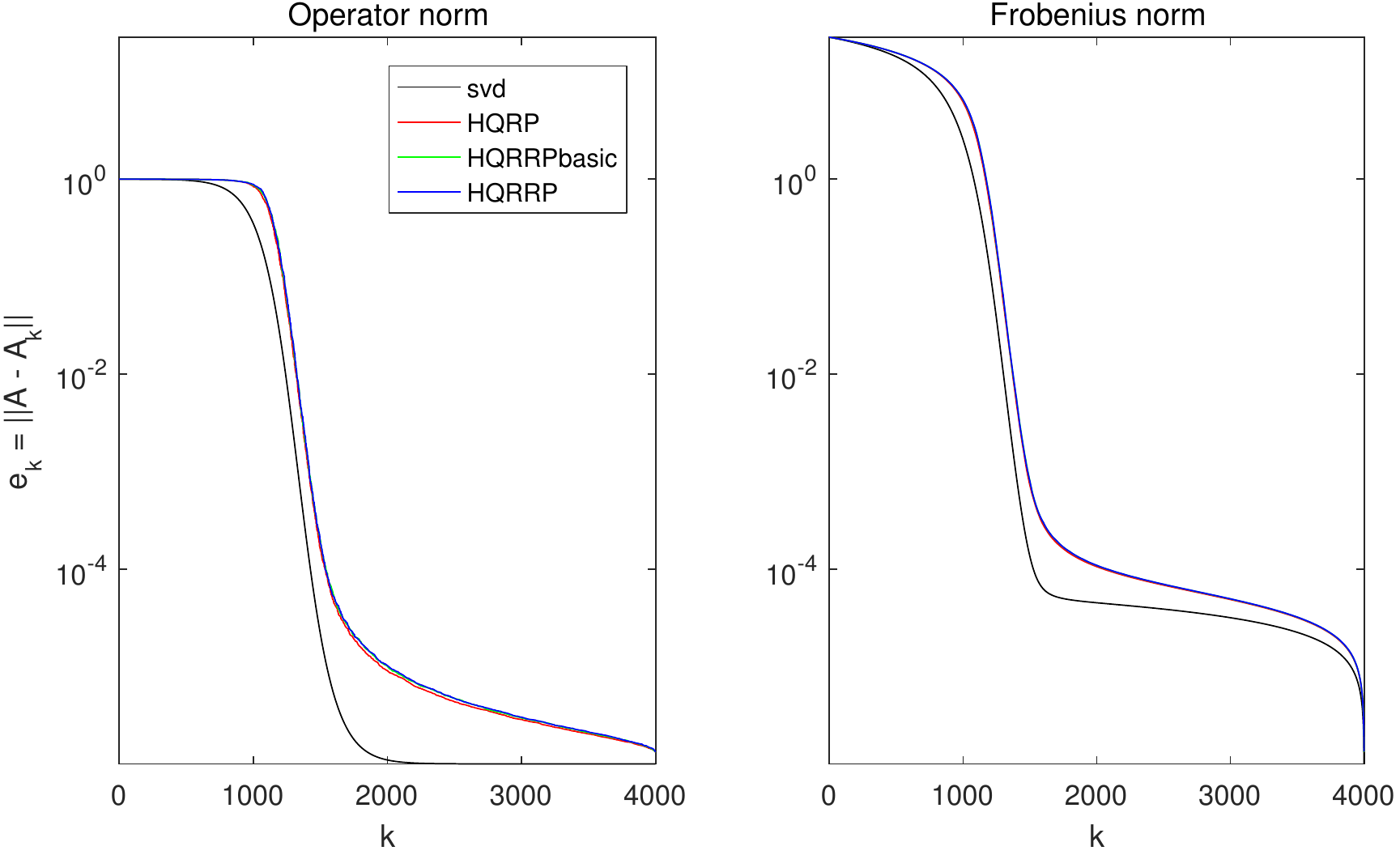}
\caption{Residual errors $e_{k}$ for ``Matrix 2'' as a function of the truncation rank $k$.
Notation is the same as in Figure \ref{fig:errors_fast}.}
\label{fig:errors_S}
\end{figure}

\begin{figure}
\includegraphics[width=\textwidth]{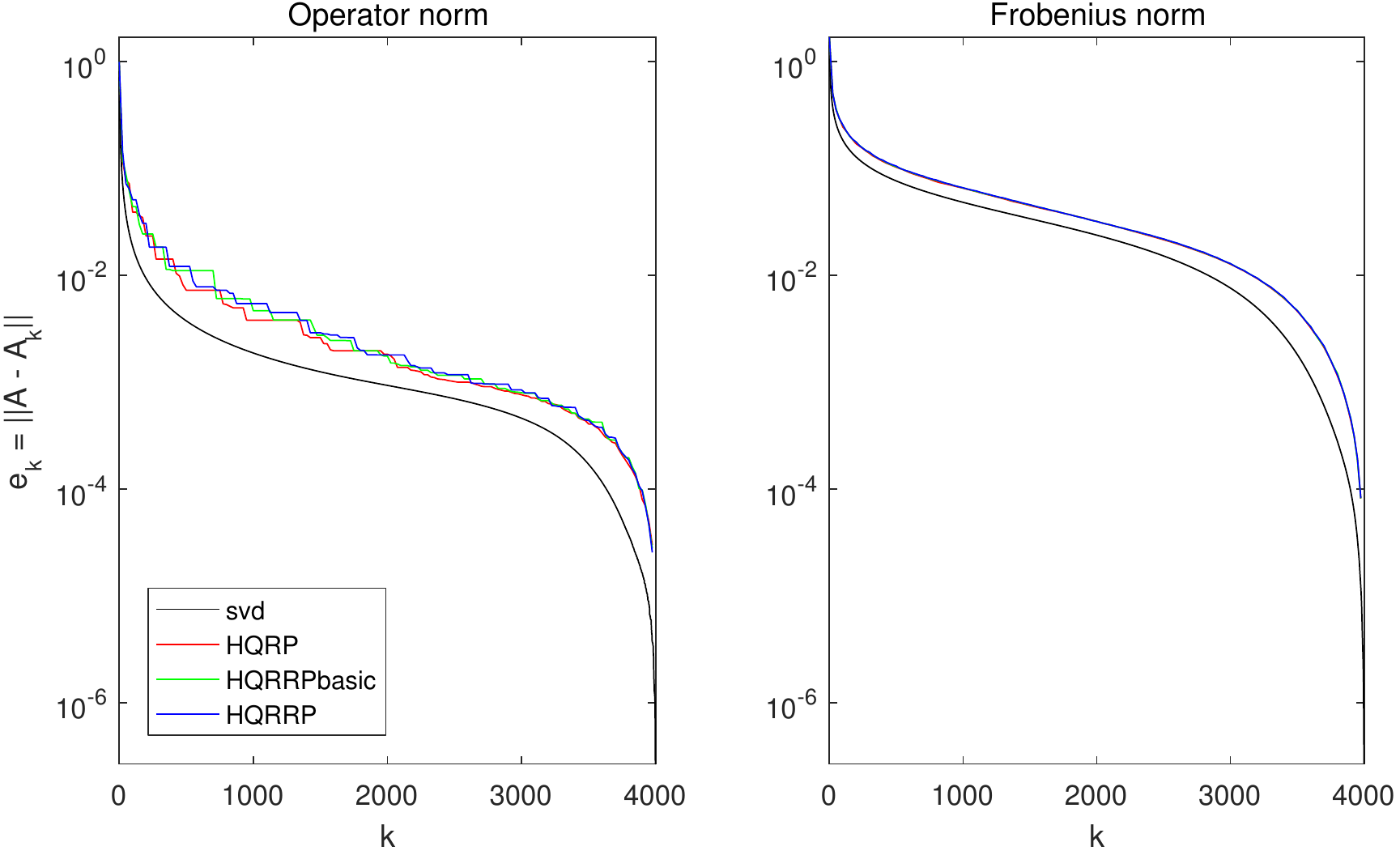}
\caption{Residual errors $e_{k}$ for ``Matrix 3'' (discretized boundary integral operator)
as a function of the truncation rank $k$. Notation is the same as in Figure \ref{fig:errors_fast}.}
\label{fig:errors_gauss}
\end{figure}

\begin{figure}
\includegraphics[width=\textwidth]{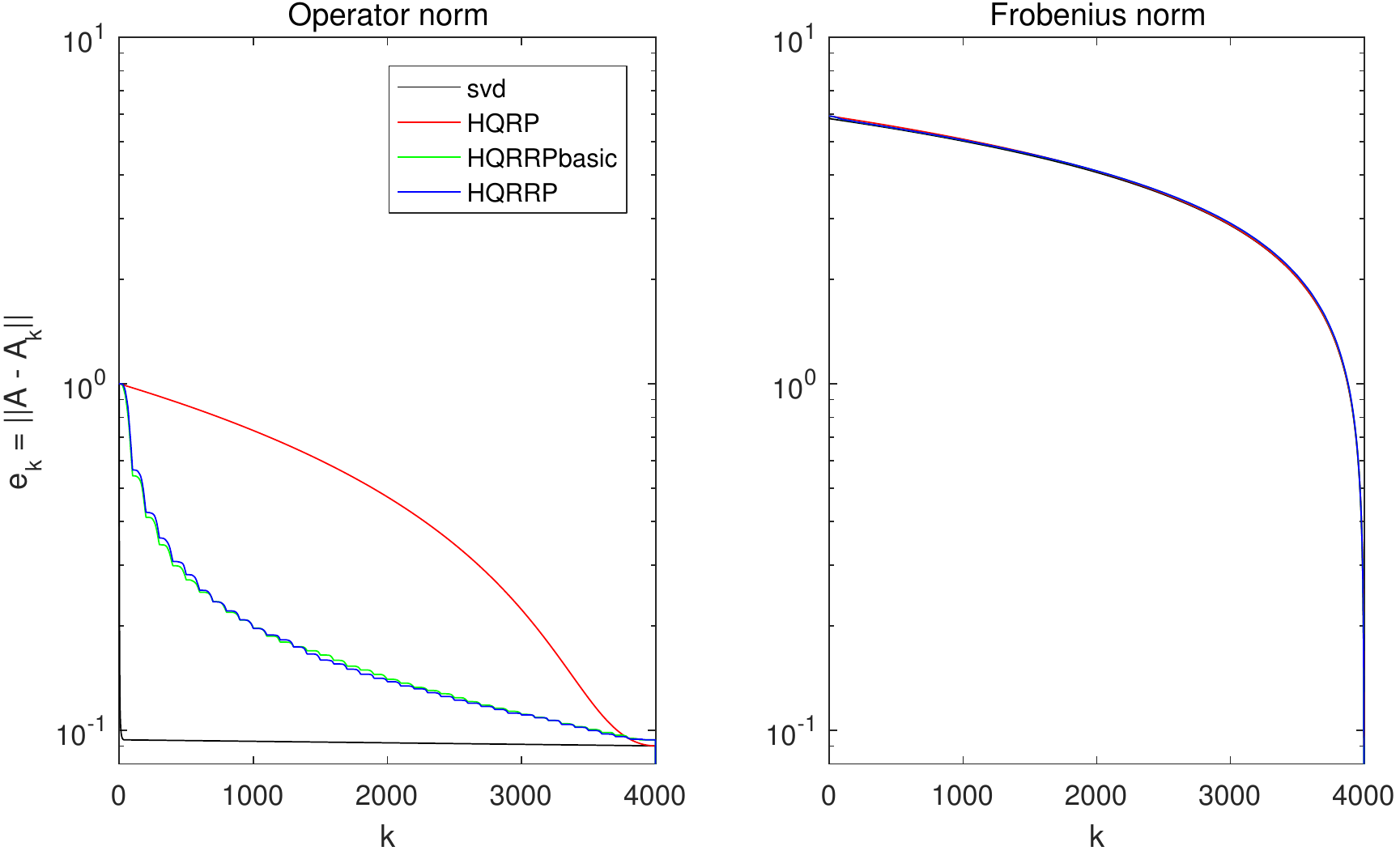}
\caption{Residual errors $e_{k}$ for ``Matrix 4'' (Kahan) as a function of the truncation rank $k$.
Notation is the same as in Figure \ref{fig:errors_fast}.}
\label{fig:errors_kahan}
\end{figure}

As a quality measure of the pivoting strategy, we computed the errors $e_{k}$
incurred when the factorization is truncated to its first $k$ components. To be
precise, these residual errors are defined via
\begin{equation}
\label{eq:def_ek}
e_{k} = \|\mtx{A}\mtx{P} - \mtx{Q}(:,1:k)\mtx{R}(1:k,:)\| = \|\mtx{R}((k+1):n,(k+1):n)\|.
\end{equation}
The results are shown in Figures \ref{fig:errors_fast} -- \ref{fig:errors_kahan}, for the
four different test matrices. The black lines in the graphs show the theoretically minimal
errors incurred by a rank-$k$ approximation. These are provided by the Eckart-Young theorem
\cite{1936_eckart_young} which states that, with $\{\sigma_{j}\}_{j=1}^{n}$
denoting the singular values of $\mtx{A}$:
\begin{center}
\begin{tabular}{ll}
$e_{k} \geq \sigma_{k+1}$& when errors are measured in the spectral norm, and\\
$e_{k} \geq \left(\sum_{j=k+1}^{n}\sigma_{j}\right)^{1/2}$& when errors are measured in the Frobenius norm.
\end{tabular}
\end{center}

We observe in all cases that the quality of the pivots chosen by the randomized method very
closely matches those resulting from classical column pivoting. The one exception is the
Kahan counter-example (``Matrix 4''), where the randomized algorithm performs much better.
(The importance of the last point should not be over-emphasized since this example is designed
specifically to be adversarial for classical column pivoting.)

\begin{figure}
\begin{tabular}{cc}
\includegraphics[width=70mm]{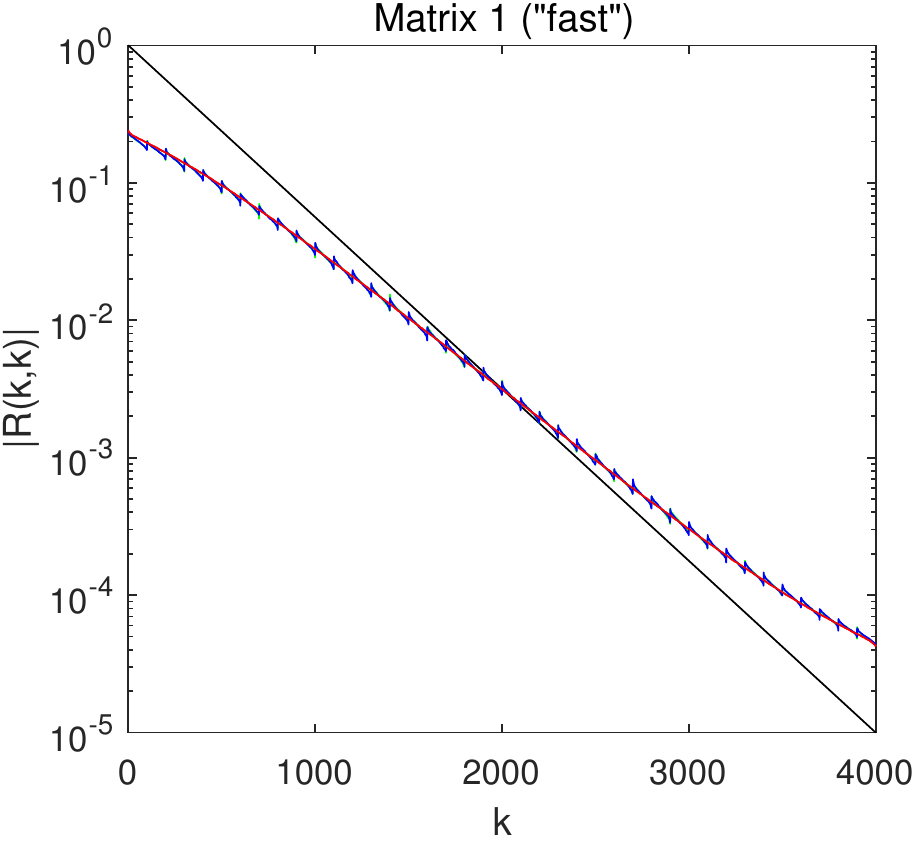} &
\includegraphics[width=70mm]{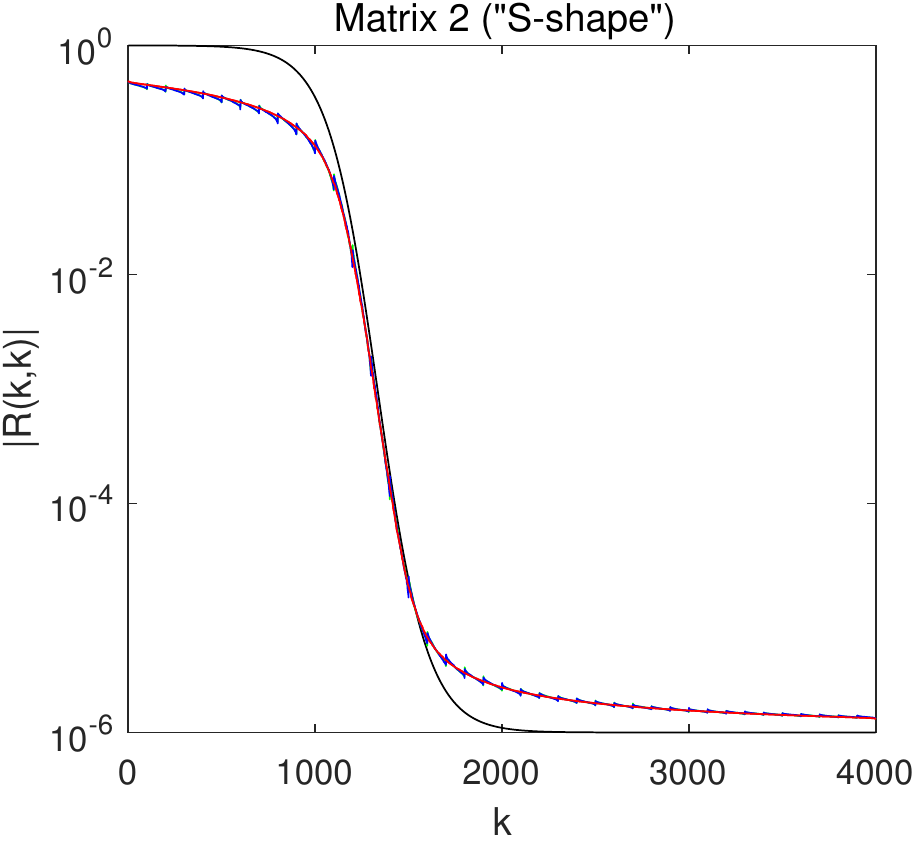} \\
\includegraphics[width=70mm]{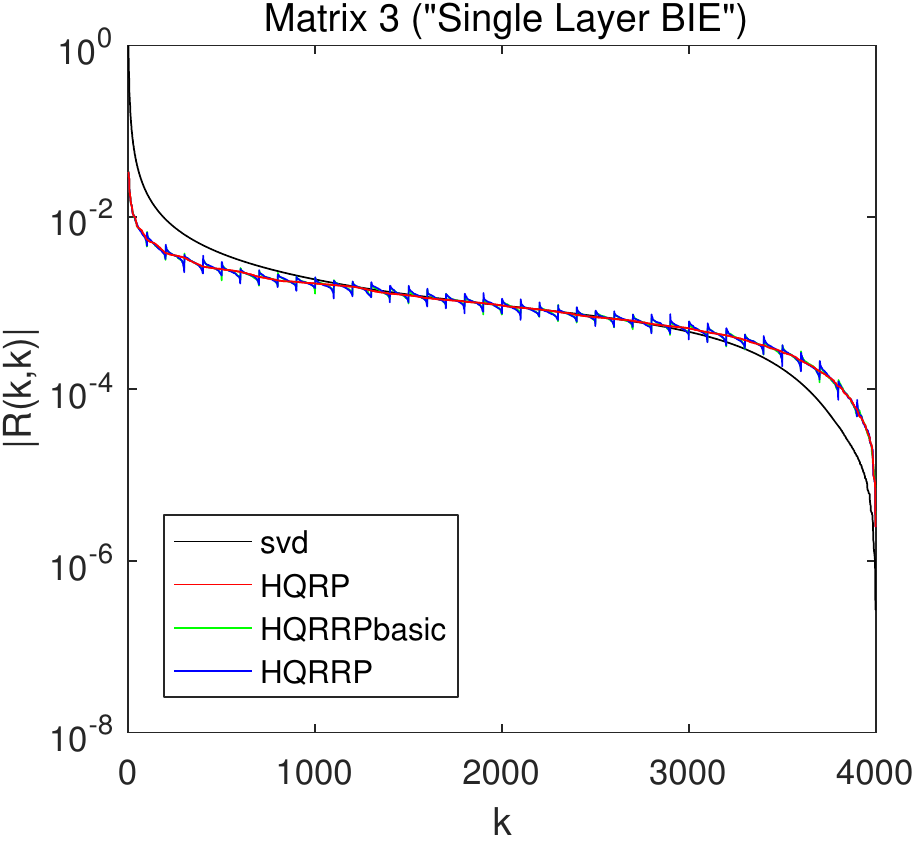} &
\includegraphics[width=70mm]{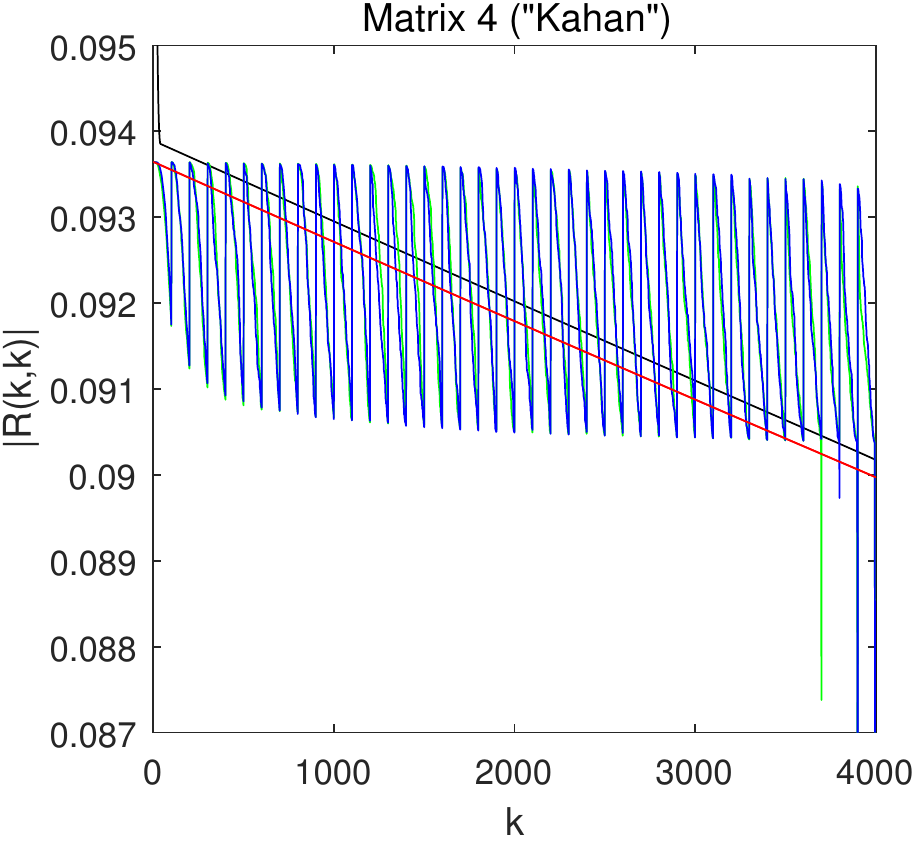}
\end{tabular}
\caption{For each of the four test matrices described in Section \ref{sec:pivotquality}, we
show the magnitudes of the diagonal entries in the ``R''-factor in a column pivoted
QR factorization. We compare classical column pivoting (red) with the two
randomized techniques proposed here (blue and green). We also show the singular values of each matrix
(black) for reference.}
\label{fig:diagonal entries}
\end{figure}

When classical column pivoting is used, the factorization (\ref{eq:AP=QRnum}) produced
has the property that the diagonal entries of $\mtx{R}$ are strictly decaying in magnitude
$$
|\mtx{R}(1,1)| \geq |\mtx{R}(2,2)| \geq |\mtx{R}(3,3)| \geq \cdots
$$
When the randomized pivoting strategies are used, this property is not enforced.
To illustrate this point, we show in Figure \ref{fig:diagonal entries} the values
of the diagonal entries obtained by the randomized strategies versus what is obtained
with classical column pivoting.

%% file: 05conclusion.tex
We have described the algorithm \HQRRP~which is a blocked version of
Householder QR with column pivoting. The main innovation compared
to earlier work is that pivots are determined in groups using a technique
based on randomized projections. We demonstrated that the quality of the chosen
pivots is for practical purposes indistinguishable from traditional
column pivoting (cf.~Figures \ref{fig:errors_fast} -- \ref{fig:errors_gauss}),
and that the dominant term in the asymptotic flop count equals that of
non-pivoted QR. Importantly, we also demonstrated through numerical
experiments that \HQRRP~is very fast, in fact almost as fast as
unpivoted \HQR.

The technique described opens up several potential avenues for
future research. The speed gains we demonstrate on single core
and shared memory multicore machines is due to the reduction in
data movement.  Equivalently, data moved between memory layers is
carefully amortized.
We expect the technique described to have an even
more pronounced advantage over traditional column pivoted QR when
implemented in more severely communication constrained environments
such as a matrix processed on a GPU or a distributed
memory parallel machine, or stored out-of-core.

The randomized sampling techniques we describe can also be used to
construct very close to optimal rank-revealing factorizations. To
describe the idea, we note that a column pivoted QR factorization
of a given matrix $A$ can be written as
\begin{equation}
\label{eq:thefuture}
A = Q\,R\,P^{*},
\end{equation}
where $Q$ is orthonormal, $R$ is upper triangular, and $P$ is a permutation
matrix. In this manuscript, we used randomized sampling to determine the
permutation matrix $P$. It turns out that for a modest additional cost,
one can build a factorization that takes the form (\ref{eq:thefuture}), but
with both $Q$ and $P$ built as products of Householder reflectors. This
generalization allows us to bring $R$ not only to upper triangular form,
but very close to being diagonal, with accurate approximations to the singular
values of $A$ on its diagonal. This discovery was reported in \cite{2015_blockQR},
and is a subject of ongoing research.

How to scale the presented algorithm to very large numbers of cores is
an open research question.  Techniques such as ``compute ahead'' will
have to be employed to ensure that the factorization of the current
panel ($ A_{11} $ and $ A_{21} $) and downdate of $ Y $ do not start
dominating the parts of the computation that can be cast in terms of \gemm. 

%% file: 06software.tex
A number of implementations of the discussed algorithm are available under 3-clause (modified) BSD license from:

\begin{center}
\texttt{https://github.com/flame/hqrrp}
\end{center}

\noindent
Included are implementations that directly link to LAPACK \cite{LAPACK3} 
as well as implementations that use the libflame \cite{inverse-siam,FLAME} library.
For those who use the LAPACK routine dgeqp3 routine, a plug compatible routine dgeqp4 is provided.

A distributed memory implementation of the algorithm has been incorporated into the Elemental
software package by Jack Poulson et al \cite{poulson2013elemental}, available at:

\begin{center}
\texttt{https://github.com/elemental/Elemental}
\end{center}

%% file: article.bbl
\begin{thebibliography}{10}

\bibitem{alpert1999hybrid}
Bradley~K Alpert.
\newblock Hybrid {G}auss-trapezoidal quadrature rules.
\newblock {\em SIAM Journal on Scientific Computing}, 20(5):1551--1584, 1999.

\bibitem{ACML}
AMD.
\newblock {AMD} {C}ore {M}ath {L}ibrary.
\newblock
  \url{http://developer.amd.com/tools-and-sdks/cpu-development/amd-core-math-library-acml/}.

\bibitem{LAPACK3}
E.~Anderson, Z.~Bai, C.~Bischof, L.~S. Blackford, J.~Demmel, Jack~J. Dongarra,
  J.~Du Croz, S.~Hammarling, A.~Greenbaum, A.~McKenney, and D.~Sorensen.
\newblock {\em LAPACK Users' guide (third ed.)}.
\newblock SIAM, Philadelphia, PA, USA, 1999.

\bibitem{RRQRblk2}
C.~H. Bischof and Gregorio Quintana-Ort\'{\i}.
\newblock Algorithm 782: Codes for rank-revealing {QR} factorizations of dense
  matrices.
\newblock {\em ACM Trans. Math. Soft.}, 24(2):254--257, 1998.

\bibitem{RRQRblk1}
C.~H. Bischof and Gregorio Quintana-Ort\'{\i}.
\newblock Computing rank-revealing {QR} factorizations of dense matrices.
\newblock {\em ACM Trans. Math. Soft.}, 24(2):226--253, 1998.

\bibitem{BiVL87}
Christian Bischof and Charles Van~Loan.
\newblock The {WY} representation for products of {Householder} matrices.
\newblock {\em SIAM J. Sci. Stat. Comput.}, 8(1):s2--s13, Jan. 1987.

\bibitem{BLAS3}
Jack~J. Dongarra, Jeremy Du~Croz, Sven Hammarling, and Iain Duff.
\newblock A set of level 3 basic linear algebra subprograms.
\newblock {\em ACM Trans. Math. Soft.}, 16(1):1--17, March 1990.

\bibitem{DDSV}
Jack~J. Dongarra, Iain~S. Duff, Danny~C. Sorensen, and Henk~A. {van der Vorst}.
\newblock {\em Solving Linear Systems on Vector and Shared Memory Computers}.
\newblock SIAM, Philadelphia, PA, 1991.

\bibitem{2006_drineas_kannan_mahoney}
Petros Drineas, Ravi Kannan, and Michael~W. Mahoney.
\newblock Fast {M}onte {C}arlo algorithms for matrices. {II}. {C}omputing a
  low-rank approximation to a matrix.
\newblock {\em SIAM J. Comput.}, 36(1):158--183 (electronic), 2006.

\bibitem{2015_blockQR_ming}
J.~Duersch and M.~Gu.
\newblock True {BLAS}-3 performance {QRCP} using random sampling, 2015.
\newblock arXiv preprint \#1509.06820.

\bibitem{1936_eckart_young}
Carl Eckart and Gale Young.
\newblock The approximation of one matrix by another of lower rank.
\newblock {\em Psychometrika}, 1(3):211--218, 1936.

\bibitem{2004_kannan_vempala}
Alan Frieze, Ravi Kannan, and Santosh Vempala.
\newblock Fast {M}onte-{C}arlo algorithms for finding low-rank approximations.
\newblock {\em J. ACM}, 51(6):1025--1041 (electronic), 2004.

\bibitem{GVL3}
Gene~H. Golub and Charles~F. Van~Loan.
\newblock {\em Matrix Computations}.
\newblock The Johns Hopkins University Press, Baltimore, 3nd edition, 1996.

\bibitem{Goto2}
Kazushige Goto and Robert~A Geijn.
\newblock Anatomy of high-performance matrix multiplication.
\newblock {\em ACM Transactions on Mathematical Software (TOMS)}, 34(3):12,
  2008.

\bibitem{Goto}
Kazushige Goto and Robert~A. van~de Geijn.
\newblock On reducing {TLB} misses in matrix multiplication.
\newblock Technical Report CS-TR-02-55, Department of Computer Sciences, The
  University of Texas at Austin, 2002.

\bibitem{gu1996}
Ming Gu and Stanley~C. Eisenstat.
\newblock Efficient algorithms for computing a strong rank-revealing {QR}
  factorization.
\newblock {\em SIAM J. Sci. Comput.}, 17(4):848--869, 1996.

\bibitem{FLAME}
John~A. Gunnels, Fred~G. Gustavson, Greg~M. Henry, and Robert~A. van~de Geijn.
\newblock {FLAME}: {F}ormal {L}inear {A}lgebra {M}ethods {E}nvironment.
\newblock {\em ACM Trans. Math. Soft.}, 27(4):422--455, December 2001.
\newblock \\Download from
  \href{http://www.cs.utexas.edu/users/flame/web/FLAMEPublications.html}{http://www.cs.utexas.edu/users/flame/web/FLAMEPublications.html}.

\bibitem{2011_martinsson_randomsurvey}
Nathan Halko, Per-Gunnar Martinsson, and Joel~A. Tropp.
\newblock Finding structure with randomness: Probabilistic algorithms for
  constructing approximate matrix decompositions.
\newblock {\em SIAM Review}, 53(2):217--288, 2011.

\bibitem{ESSL}
IBM.
\newblock {E}ngineering and {S}cientific {S}ubroutine {L}ibrary.
\newblock \url{http://www-03.ibm.com/systems/power/software/essl/}.

\bibitem{MKL}
Intel.
\newblock {M}ath {K}ernel {L}ibrary.
\newblock \url{http://developer.intel.com/software/products/mkl/}.

\bibitem{Joffrain:2006:AHT:1141885.1141886}
Thierry Joffrain, Tze~Meng Low, Enrique~S. Quintana-Ort\'{\i}, Robert
  {v}an~{d}e {G}eijn, and Field~G. {V}an {Z}ee.
\newblock Accumulating {H}ouseholder transformations, revisited.
\newblock {\em ACM Trans. Math. Softw.}, 32(2):169--179, June 2006.

\bibitem{1966_kahan_NLA}
William Kahan.
\newblock Numerical linear algebra.
\newblock {\em Canadian Math. Bull}, 9(6):757--801, 1966.

\bibitem{2007_martinsson_PNAS}
Edo Liberty, Franco Woolfe, Per-Gunnar Martinsson, Vladimir Rokhlin, and Mark
  Tygert.
\newblock Randomized algorithms for the low-rank approximation of matrices.
\newblock {\em Proc. Natl. Acad. Sci. USA}, 104(51):20167--20172, 2007.

\bibitem{mahoney2011randomized}
Michael~W Mahoney.
\newblock Randomized algorithms for matrices and data.
\newblock {\em Foundations and Trends{\textregistered} in Machine Learning},
  3(2):123--224, 2011.

\bibitem{2015arXiv150307157M}
P.-G. {Martinsson} and S.~{Voronin}.
\newblock A randomized blocked algorithm for efficiently computing
  rank-revealing factorizations of matrices, mar 2015.
\newblock ArXiv.org e-print \#1503.07157. To appear in the SIAM Journal on
  Scientific Computing.

\bibitem{2006_martinsson_random1_orig}
Per-Gunnar Martinsson, Vladimir Rokhlin, and Mark Tygert.
\newblock A randomized algorithm for the approximation of matrices.
\newblock Technical Report Yale CS research report YALEU/DCS/RR-1361, Yale
  University, Computer Science Department, 2006.

\bibitem{2011_martinsson_random1}
Per-Gunnar Martinsson, Vladimir Rokhlin, and Mark Tygert.
\newblock A randomized algorithm for the decomposition of matrices.
\newblock {\em Appl. Comput. Harmon. Anal.}, 30(1):47--68, 2011.

\bibitem{2015_blockQR}
P.G. Martinsson.
\newblock Blocked rank-revealing qr factorizations: How randomized sampling can
  be used to avoid single-vector pivoting, 2015.
\newblock arXiv preprint \#1505.08115.

\bibitem{poulson2013elemental}
Jack Poulson, Bryan Marker, Robert~A Van~de Geijn, Jeff~R Hammond, and
  Nichols~A Romero.
\newblock Elemental: A new framework for distributed memory dense matrix
  computations.
\newblock {\em ACM Transactions on Mathematical Software (TOMS)}, 39(2):13,
  2013.

\bibitem{inverse-siam}
Enrique~S. Quintana, Gregorio Quintana, Xiaobai Sun, and Robert van~de Geijn.
\newblock A note on parallel matrix inversion.
\newblock {\em SIAM J. Sci. Comput.}, 22(5):1762--1771, 2001.

\bibitem{QRP:SIAM}
Gregorio Quintana-Ortí, Xiaobai Sun, and Christian~H. Bischof.
\newblock A {BLAS}-3 version of the {QR} factorization with column pivoting.
\newblock {\em SIAM Journal on Scientific Computing}, 19(5):1486--1494, 1998.

\bibitem{ScVL89}
Robert Schreiber and Charles Van~Loan.
\newblock A storage-efficient {WY} representation for products of {H}ouseholder
  transformations.
\newblock {\em SIAM J. Sci. Stat. Comput.}, 10(1):53--57, Jan. 1989.

\bibitem{1998_stewart_volume1}
G.W. Stewart.
\newblock {\em Matrix Algorithms Volume 1: Basic Decompositions}.
\newblock SIAM, 1998.

\bibitem{libflame_ref}
Field~G. {V}an {Z}ee.
\newblock {\em {\tt libflame}: {T}he {C}omplete {R}eference}.
\newblock {\tt www.lulu.com}, 2012.
\newblock \\Download from
  \href{http://www.cs.utexas.edu/users/flame/web/FLAMEPublications.html}{http://www.cs.utexas.edu/users/flame/web/FLAMEPublications.html}.

\bibitem{CiSE09}
Field~G. {V}an {Z}ee, Ernie Chan, Robert van~de Geijn, Enrique~S.
  Quintana-Ort\'{\i}, and Gregorio Quintana-Ort\'{\i}.
\newblock The libflame library for dense matrix computations.
\newblock {\em IEEE Computation in Science \& Engineering}, 11(6):56--62, 2009.

\bibitem{BLIS1}
Field~G. {Van Zee} and Robert~A. {v}an~{d}e {G}eijn.
\newblock {BLIS}: A framework for rapidly instantiating {BLAS} functionality.
\newblock {\em {ACM} Transactions on Mathematical Software}, 41(3), 2015.

\bibitem{2014_martinsson_CUR}
S.~Voronin and P.G. Martinsson.
\newblock A {CUR} factorization algorithm based on the interpolative
  decomposition.
\newblock {\em ar{X}iv.org}, 1412.8447, 2014.

\bibitem{woodruff2014sketching}
David~P Woodruff.
\newblock Sketching as a tool for numerical linear algebra.
\newblock {\em arXiv preprint arXiv:1411.4357}, 2014.

\bibitem{xianyi2012model}
Xianyi Zhang, Qian Wang, and Yunquan Zhang.
\newblock Model-driven level 3 {BLAS} performance optimization on loongson 3a
  processor.
\newblock In {\em Proceedings of the 2012 IEEE 18th International Conference on
  Parallel and Distributed Systems}, pages 684--691. IEEE Computer Society,
  2012.

\end{thebibliography}
